\documentclass{article}

\usepackage[T1]{fontenc}
\usepackage[utf8]{inputenc}
\usepackage[cm]{aeguill}
\usepackage[english]{babel}

\usepackage{amsmath}
\usepackage{amssymb}
\usepackage{amsthm}
\usepackage{nicefrac}
\usepackage{a4wide}

\usepackage[dvips]{graphicx} %changement important ?

\def\proc#1{\medbreak\noindent{\it #1}\hspace{1ex}\ignorespaces}
\def\ep{\noindent{\hfill $\Box$}}

\newcommand{\eps}{\varepsilon}
\renewcommand{\phi}{\varphi}

\DeclareMathOperator{\Card}{Card}
\DeclareMathOperator{\Stab}{Stab}
\DeclareMathOperator{\rk}{rk}
\DeclareMathOperator{\vect}{Vect}
\DeclareMathOperator{\mat}{Mat}

\newcommand{\ie}{{\emph{i.e.~}}}

\newcommand{\R}{\mathbb{R}}
\newcommand{\Q}{\mathbb{Q}}
\newcommand{\N}{\mathbb{N}}
\newcommand{\Z}{\mathbb{Z}}

\newcommand{\sg}{\mathrm{Sing}}
\newcommand{\reg}{\mathrm{Reg}}
\newcommand{\Pron}{\mathcal{P}}
\newcommand{\tP}{\widetilde{\mathcal{P}}}
\newcommand{\tA}{\widetilde{A}}

\newcommand{\tq}{\,;\,}

\newcommand{\nr}[1]{\left\Vert #1\right\Vert}
\newcommand{\abs}[1]{\left\vert #1\right\vert}
\newcommand{\restr}[2]{\left.#1\right|_{#2}}

\newcommand{\Aron}{\mathcal{A}}
\newcommand{\Bron}{\mathcal{B}}
\newcommand{\Cron}{\mathcal{C}}
\newcommand{\Dron}{\mathcal{D}}

\newtheorem{theo}{Theorem}[section]
\newtheorem{prop}[theo]{Proposition}
\newtheorem{lemma}[theo]{Lemma}
\newtheorem{cor}[theo]{Corrolary}

\newenvironment{definition}{\refstepcounter{theo}
\proc{Definition~\thetheo. }}{\medbreak}

\newenvironment{rem}{\refstepcounter{theo}
\proc{Remark~\thetheo.}}{\medbreak}

\title{Complexity and cohomology for cut and projection tilings}

\author{Antoine \sc{Julien}} %\affil{1}}

\date{
\small\textit{Universit\'e de Lyon; Universit\'e Claude Bernard Lyon 1;} \\
\small\textit{Institut Camille Jordan CNRS UMR 5208;}\\
\small\textit{B\^atiment Braconnier}\\
\small\textit{43 boulevard du 11 novembre 1918, F-69622 Villeurbanne Cedex.} \\
\small\textit{e-mail:} \texttt{julien@math.univ-lyon1.fr} \medskip \\
April 2008}

\begin{document}

\maketitle

\begin{abstract}
We consider a subclass of tilings, the tilings obtained by cut and projection.
Under somewhat standard assumptions, we show that the natural complexity
function has polynomial growth. We compute its exponent $\alpha$ in terms of the
ranks of certain groups which appear in the construction. We give bounds for
$\alpha$.
These computations apply to some well known tilings, such as the octagonal
tilings, or tilings associated with billiard sequences.
A link is made between the exponent of the complexity, and the fact
that the cohomology of the associated tiling space is finitely generated
over $\Q$.
We show that such a link cannot be established for more general tilings, and
we present a counter-example in dimension one.
\end{abstract}

\section*{Introduction}

A major motivation for studying aperiodic tilings came from the
physics of quasi-crystals. Even though the aperiodic Wang tilings and Penrose
tilings have already been studied in the 1960's and the 1970's, respectively,
most of the work on this subject was made since the first discovery of a
quasi-crystalline material in the 1980's.
Such a material appeared to be very ordered, because its diffraction pattern
showed clear peaks. However, it could not be periodic, as the diffraction
pattern had symmetries which should not have occurred according to the
classical periodic models of crystals.
Until this observation, a crystalline material was seen as a set of atoms, the
position of which was determined by a regular lattice. The group of
transformations of $\R^3$ preserving the lattice determined the symmetries of
the crystal.
The appearance of quasi-crystals in physics created the need
for new mathematical objects in order to model these materials.

The cut and projection method plays a major role for the description of such
aperiodic yet ordered sets. The general idea is to project a ``slice'' of a
higher-dimensional regular lattice of $\R^N$ on a subspace of dimension $d$.
The resulting point pattern should not be periodic, but should inherit the
symmetry properties of the bigger lattice.
This method actually provides important examples of tilings, such as the
octagonal tilings, the dodecagonal tilings, or the icosahedral tilings.

The notion of ordered point set, or ordered tiling, is ambiguous. The
definition itself of what should be an ordered tiling varies accordingly to the
author's preoccupations.
A most common assumption is the notion of finite local complexity (FLC). One
requires that there are only finitely many local configurations in the tiling
(up to translation). This is a qualitative condition on complexity.
A dynamical definition of order could be recurrence of orbits for the associated
dynamical system. In combinatorial terms, this corresponds to a property of
repetitivity: every finite patch of the tiling appears infinitely
often, and within a prescribed range depending on the size of the patch.
Another dynamical condition would be to require that the entropy of the
dynamical system is zero. In terms of ``patch-counting entropy'', it is
equivalent to require that the growth of the patch-counting function (which to
$R$ associates the number of patches included in a ball of size $R$) is
sub-exponential.

Our approach to understand the notion of ordered point pattern involves a
complexity function, as in the entropy case, but with greater precision.
The notion of complexity is well defined for bi-infinite words over a finite
alphabet: the function $n \mapsto p(n)$ associated to a word $w$ counts the
number of subwords of $w$ of length $n$.
This definition can be extended to multi-dimensional tilings, though not in a
canonical way: $p(n)$ should count (up to translation) the number of patches
of size $n$. This definition depends on a suitable definition of ``patch of size
$n$''.
While many properties are known about the complexity of one-dimensional words
(that is for one-dimensional tilings), the situation is not as well known
for multi-dimensional tilings.

In this article, we define a notion of ``patch of size $n$'' in the
framework of cut and projection tilings, in order to define a complexity
function. We compute this function for cut and projection tilings with
almost-canonical acceptance domain. Almost-canonical acceptance domains, as we
will define, are a generalization of the canonical case, when the
projected ``slice'' of the lattice is a cylinder, the base of which is the
projection of the unit cube.
Almost canonical acceptance domains are domains for which it is possible to
define cut hyperplanes in the ``internal space'' (a supplementary of the
``physical space'').

We will prove (theorem~\ref{theorem}) that for a cut and projection tiling of
dimension $d$ arising as the projection of the regular lattice $\Z^N$ with
almost-canonical acceptance domain, the complexity $p(n)$ grows like $n^\alpha$,
where $\alpha$ is an integer which depends on the data of the cut and projection
method.
Namely, it depends on the rank of the stabilizers of the cut planes,
and will be computed in examples.
We give bounds for $\alpha$: for a tiling with no period, we have
$d \leq \alpha \leq d(N-d)$ (theorems~\ref{theorem} and~\ref{theorem-cohom}).
Generically, the complexity is maximal, that is $\alpha = d (N-d)$.

Another notion of ``ordered structure'' could be understood in terms of the
topology of the tiling space.
Given a tiling and its associated tiling space $\Omega$, one can compute its
\v{C}ech cohomology.
A qualitative distinction can be done between tiling spaces with
finitely generated rational cohomology, and tiling spaces with infinitely
generated rational cohomology.
This distinction seems to characterize the difference between a topologically
complex space and a topologically simpler one, and our intuition is that
the class of tilings with finitely generated rational cohomology should have
some specific interest, being a class of more ``ordered'' tilings. Note that the
cohomology groups of substitution tiling spaces are finitely generated over
$\Q$ (see~\cite{AP}), as well as all the tilings which have been used to model
quasi-crystals in physics.

For cut and projection tilings, there is actually a link between these two
notions of order: the cohomological ``complexity'' and the growth of the
complexity function $p$.
We prove in section~\ref{sec:links} (theorem~\ref{theorem-cohom}) that, in the
case of a $d$-dimensional cut and projection tiling with almost-canonical
acceptance domain and with no periods, the following equivalence holds: the
cohomology of the tiling space is finitely generated if and only if the
complexity of the tiling is low, that is $p(n)$ grows like $n^d$.

Such a link between high complexity and infinitely generated cohomology does
not hold in a more general framework.
We prove that if a word is not too complex (which means that its complexity is
dominated by $C n$ with $C$ a constant), then its associated tiling space has
finitely generated cohomology.
This result is in line with the statement above, but the converse does not hold
in full generality: we define a word, the complexity of which is exponential,
but such that its associated tiling space has finitely generated cohomology.
In order to prove these results, we prove that the tiling space of a word is
homeomorphic to the inverse limit of its Rauzy graphs.
It was already known that the space was an inverse limit of graphs, thanks to
the Anderson-Putnam-G\"ahler complexes (see~\cite{Sad}). However, Rauzy
graphs are more suited than G\"ahler complexes to understand the combinatorics
of a word. Our proof uses the now classical methods developed in~\cite{AP},
or~\cite{Sad}.

The questions raised in this article have interested mathematicians for some
time. Answers for some of these problems are already known in specific cases.
The case of complexity for one-dimensional cut and projection tilings,
also known as cubic billiard sequences was already investigated by~\cite{AMST}
in dimension 3, and then by~\cite{Bar} in any dimension. In these articles,
an exact and explicit polynomial formula is given for the complexity function,
assuming certain additional arithmetic condition.
It seems however, that the complexity problem was not yet solved in a general
framework.
Berth\'e and Vuillon obtained more precise results than ours for the complexity,
in specific case of $(3,2)$, $(N,N-1)$ and to a certain extend $(N,2)$-cut
and projection tilings (or ``discrete planes''), see~\cite{BV} and~\cite{Vui}.
Their method consists in modifying the tiling in order to obtain a
$\Z^2$-subshift, and study its rectangular complexity.
The method used in our paper seems very similar to the one used in~\cite{AMST}
in the sense that the problem of counting patches (or words in their case) is
reduced to the
problem of counting connected components of a certain set (counting cells of a
certain decomposition of the two torus, respectively).

There are also other lines of approach to study complexity in more
than one dimension: Lagarias and Pleasants~(\cite{LP}) have worked on the
complexity of general Delone sets of $\R^d$. They prove results of minimal
complexity. More precisely, they state the following result: assume the
complexity function $p$ of a Delone set $\Dron$ satisfies $p(n) \leq C n$ for
some $n \in \N$ and for some constant $C$ depending (in a weak sense) on the
set. Then $\Dron$ is an ideal crystal, that is, its group of periods has rank
$d$.
We prove here a result which is related to a conjecture in their
article, even though we prove it for the much more specific case of cut and
projection Delone sets: the group of periods of a cut and projection Delone
set of $\R^d$ is of rank at most $k$, if and only if $p(n) / n^{d-k}$ is
bounded from below by a positive constant.
However, we do not have control on the constant involved, and so this does not
yet fully confirm the conjecture for almost-canonical cut and projection
tilings.

\section{Definitions and General Properties}

\subsection{Cut and Projection Tiling Spaces}

Let $\R^N$ be the standard Euclidean space with basis $(e_i)_{i=1}^N$, and
distinguished sublattice $\Z^N$.
In the following, this space will be normed by the $1$-norm, which we
will note $\nr{.}_1$, and the closed balls will be noted
$\Bron(x,r):=\{y \in \R^N \tq \sum_i{\abs{x_i - y_i}} \leq r\}$, for $r > 0$
and $x \in \R^N$.

Let $E$ be a linear subspace of $\R^N$, $F$ a complementary subspace, $\pi$
and $\pi_F$ be the projections respectively on $E$ and $F$, associated to the
decomposition $\R^N = E \oplus F$.
We note $d$ the dimension of $E$, and we make the following standard assumption:
\proc{Assumption.}
\label{hyp}
The restriction $\restr{\pi}{\Z^N}$ is one-to-one, and $\pi_F \big(\Z^N \big)$
is a dense subgroup of $F$.
\medbreak

\begin{definition}
We define the following objects:
\begin{itemize}
 \item The \emph{acceptance domain} $K$ is a compact subset of $F$, which is
the closure of its interior, and which we assume to be a convex polytope (the
closed convex hull of a finite set).
 \item $V$ the (finite) set of vertices of $K$.
 \item The \emph{acceptance strip}, which we note $S$, is the set $K+E$.
 \item $\Gamma := \pi_F (\Z^N)$ the projection of $\Z^N$ on $F$, which we
assumed to be a dense subgroup of $F$ (by assumption~\ref{hyp}).
\end{itemize}
\end{definition}

Let $C$ be the closed unit cube, \ie the convex set $C = \{\sum{\lambda_i e_i}
\tq 1 \leq i \leq N, \, 0 \leq \lambda_i \leq 1 \}$.
In the case $K = \pi_F (C)$, the cut and projection method is said to have
\emph{canonical acceptance domain}.
One can read the article~\cite{ODK} for general results about cut and
projection tilings with canonical acceptance domain, and their relation
with other construction methods.

In a more general setting, we assume that the acceptance domain is
\emph{almost-canonical} in the following sense. From now on, we will assume
that the acceptance domain is almost canonical.
\begin{definition}\label{hypK}
The acceptance domain is said to be \emph{almost-canonical} if the following
holds:
$K$ is a convex polytope, such that $(\partial K + \Gamma)$ is a countable
union of hyperplanes of $F$.
In addition, given any $(N-d-1)$-dimensional face $f$ of $K$, and any vertex
$v$ of $f$, and $H_f$ the linear hyperplane parallel to $f$, we assume that
there is a neighborhood of $0$ in $H_f$ which can be covered by finitely many
translates of $f$ of the form $(f-v) + \gamma$, with
$\gamma \in \Gamma \cap H_f$.
\end{definition}

Remark that $f-x$ has non empty interior in $H_f$. Therefore, if $H_f \cap
\Gamma$ is dense in $H_f$ for all hyperplanes $H_f$, then the acceptance domain
$K$ is almost-canonical.
We will prove later that a canonical acceptance domain is almost-canonical. An
almost-canonical domain is a generalization of a canonical domain which allows
to define cut planes, and such that the singular points (cf. infra) have the
same structure as in the canonical case.

\begin{definition}
The set $(\R^N = E \oplus F, \Z^N, \pi, K)$ is called the \emph{data} of a
\emph{$(N,d)$-cut and projection method}.
\end{definition}

\begin{rem}
Sometimes, the space $E$ (the one which will be tiled), is referred to as the
\emph{physical space}, whereas the space $F$ (which can be seen as a
parameter space), is referred to as the \emph{internal space}.
\end{rem}

A cut and projection method can allow one to construct tilings. Namely, in the
case of a canonical acceptance domain, a very elegant and canonical way to do
so is described in~\cite{ODK}.
However, we will not describe cut and projection \emph{tilings} but rather cut
and projection \emph{Delone sets}, which amounts to the same when it comes to
computing the complexity.

A Delone set of $\R^d$ is a set $\Pron$ which is uniformly discrete and
relatively dense.
Uniform discreteness means that there exists $r > 0$ such that for all
$x \in \Pron$, $\Bron(x,r) \cap \Pron = \{x\}$.
Relative density means that there exists $R >0$ such that for all $x \in
\R^d$, the set $\Bron(x,R) \cap \Pron$ is not empty.
Such a set is called a $(r,R)$-Delone set.

\begin{figure}[htp]
\begin{center}
\includegraphics{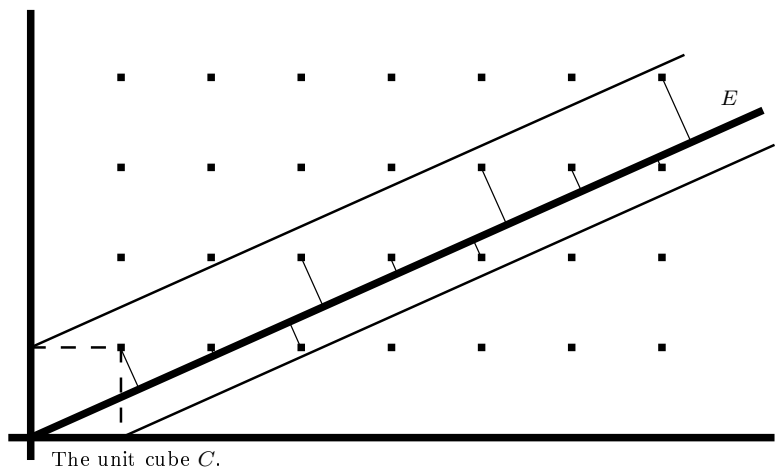}
\caption{A one-dimensional cut and projection tiling: elements of the
``slice'' $C+E$ are projected on $E$.}
\label{fig-cp}
\end{center}
\end{figure}

We will construct Delone sets by projection on $E$ of a subset of the
lattice $\Z^N$ determined by the acceptance domain $K$. The general idea is
to project the ``slice'' $S \cap \Z^N$ on $E$, as in figure~\ref{fig-cp}, but we
have to be a little bit careful in the case $\Z^N$ intersects the border of $S$.

\begin{definition}
We say that $x \in \R^N$ is a \emph{singular} point ($x \in \sg$), if
$(x + \Z^N) \cap \partial S \neq \varnothing$, where $\partial S$ stands for
the boundary of the acceptance strip $S$ in $E$.
If $x \notin \sg$, we say that $x$ is \emph{regular} and we note $x \in \reg$.
\\
Given $n \in \N \cup \{0\}$, we say that $x \in S$ is $n$-singular if
$(x + \Z^N) \cap \Bron (x,n) \cap \partial S \neq \varnothing$. We then note
$x \in \sg_n$.
We say that $x \in S$ is $n$-regular if it is not $n$-singular, and we note
$x \in \reg_n$.
\end{definition}

When a point $x$ is singular, there is some $y \in (\Z^N + x) \cap
\partial S$, so that it is unclear whether $y$ should be considered as inside
or outside $S$.
This could cause an ambiguity in the definition of the cut and projection
Delone set associated to $x$. This is why in the definition below, we only
consider Delone set associated to regular points.
The $n$-regular points are an adaptation of the regular points: they are
points for which it is possible to construct without ambiguity a patch of size
$n$ in the following sense.
\begin{definition}
For all $x \in \reg$, we define the following objects:
\begin{itemize}
\item The lifted set associated to $x$ is:
$\tP (x) := (x + \Z^N) \cap S$;
\item The Delone set associated to $x$ is: $\Pron (x) := \pi \big( \tP(x) \big)
\subseteq E$.
\end{itemize}
For all $n \in \N \cup \{0\}$, for all $x \in \reg_n$, we define:
\begin{itemize}
\item The lifted patch of size $n$ associated to $x$ (or of \emph{seed} $x$):
$\tA(x,n) = (x+\Z^N) \cap \Bron(x,n) \cap S$;
\item The patch of size $n$ associated to $x$:
$A(x,n) := \pi\big(\tA(x,n)\big)$.
\end{itemize}
Finally, we define the pointed lifted patch $\tA_{\mathrm{pt}}(x,n)$ as
$\tA(x,n)$ with distinguished point $x$.
Similarly, the pointed patch $A_{\mathrm{pt}}(x,n)$ is $A(x,n)$ with
distinguished point $\pi(x)$.
\end{definition}

The $\Pron (x)$ are Delone sets of $E$ which are known to be aperiodic and
repetitive. More precisely, under our assumptions, the group of periods of
$\Pron (x)$ is never a regular lattice of $E$.
If $E \cap \Z^N = \{0\}$, the group of periods is trivial, and the resulting
Delone sets are said to be completely aperiodic.

The $A(x,n)$ are bounded sets.
One can be tempted  to call $A(x,n)$ the ``patch of center $x$ and radius
$n$'', but this could be misleading: $x \in \R^N$ whereas
$A(x,n) \subset E$. We call $x$ a \emph{seed} of $A(x,n)$.
Actually, the ``physical'' center of $A(x,n)$ in $E$ is $\pi (x)$: in the
decomposition $x = \pi (x) + \pi_F (x)$, the term $\pi (x)$
encodes the position of $A(x,n)$ in $E$, and $\pi_F (x)$ encodes the shape
of the patch. The question ``how does $\pi_F (x)$ encode the shape of
$A(x,n)$'' is all this paper is about.

Remark that we chose $x \in S$ for the definition of $n$-singular points, and
hence for the definition of patches of size $n$. This is because if $n$ is
given, and one takes $x$ far away from $E$, one can have that $A(x,n)$ is
empty. As we want to consider patches of size $n$ and not patches of size
at most $n$, we take $x \in S$ in the definition of the patches.

If $x \in \reg$, the patch $A (x,n)$ is a subpatch of $\Pron (x)$,
as it can be extended: it is actually a subset.
If $x \in \big( \reg_n \setminus \reg \big)$, $A (x,n)$ is a subpatch of
$\Pron (y)$ for all $y$ (by repetitivity), that is \emph{up to translation} by
an element of $E$, $A(x,n)$ appears in $\Pron(y)$.

\begin{rem}\label{remark-centers}
There is \emph{a priori} a difference between patches and pointed patches.
Beware that a patch of radius $n$ could have several centers:
there could exist $x,y$ such that $\pi(x) \neq \pi(y)$, and
$A(x,n) = A(y,n)$.
That is, a patch could correspond to several pointed patches.
However, the number of candidates to be the center of any patch of the form
$A(x,n)$ is bounded, and this bound only depends on the data of the cut and
projection method.
This fact will be proved in lemma~\ref{lemma-center}.
\end{rem}

We can now define the complexity function, the function which counts the number
of different patches of a given size:
\begin{definition}
Let $\Pron (x)$ be a Delone set. The complexity function is the function
$p : \N \rightarrow \N$, such that $p(n)$ is the number of subpatches of size
$n$ of $\Pron (x)$, up to translation by an element of $E$.
We can write it:
\[
p(n) := \Card \Big( \{ A(y,n) \tq y \in \tP(x) \} / E \Big)
\]
\end{definition}

The repetitivity ensures the following property:
\begin{prop}
The complexity $p$ of a given Delone set $\Pron (x)$ equals the complexity of
any Delone set $\Pron (x')$ obtained with the same cut and projection data.
%It equals the complexity of the subshift (that is of all $A(x,k)$ for
%$x \in \reg_k$, and not only $x \in \reg$).
We deduce:
\[
p(n) = \Card \Big( \{ A(y,n) \tq y \in \reg_n \} / E \Big)
\]
Furthermore, as we count the patches only up to translation by an element of
$E$, we have:
\[
p(n) = \Card \Big( \{ A(y,n) \tq y \in K \cap \reg_n \} / E \Big)
\]
\end{prop}

Similarly, we define the pointed complexity function, which counts the number
of pointed patches up to a translation of $E$:
\begin{definition}
The pointed complexity function associated to the cut and projection method is:
\[
p_\mathrm{pt}(n) = \Card \Big( \{ A_\mathrm{pt}(y,n)
\tq y \in K \cap \reg_n \} \Big)
\]
\end{definition}

The complexity function and pointed complexity function are linked by the
following inequalities:
\begin{prop}
\label{prop-center}
There exists $\lambda > 0$, which only depends on the data of the cut and
projection method, such that for all $n \in \N$:
\[
\lambda p_\mathrm{pt}(n) \leq p(n) \leq p_\mathrm{pt}(n)
\]
\end{prop}

We first state the following lemma, which proves the assertion of
remark~\ref{remark-centers}.
\begin{lemma}
\label{lemma-center}
Given a cut and projection method, there exists $M \in \N$ such that for all
$n \in \N$, for all $x \in \reg_n$, the patch $A(x,n)$ has at most $M$ centers;
that is:
\[
\Card{\Big( \big\{ y \in \tA(x,n) \tq A(x,n) = A(y,n) \big\}\Big)} \leq M
\]
\end{lemma}

\proc{Demonstration.}
The size of a patch is given by the norm in $\R^N$. We want to express the size
of a patch of size $n$ in terms of the norm induced in $E$ by the norm of
$\R^N$.
Remark that if $n \in \N$ and $x \in K$, then for all $y \in \tP(x)$,
\[
\nr{\pi(y) - \pi(x)}_1 \leq \nr{\pi(y) - y}_1 + \nr{y-x}_1 + \nr{x-\pi(x)}_1
\]
and
\[
\nr{y-x}_1 \leq \nr{\pi(y)-y}_1 + \nr{\pi(y)-\pi(x)}_1 + \nr{x-\pi(x)}_1
\]
But $\pi(y) - y \in K$ and $K$ is compact, so there exists an upper bound to
$\nr{\pi(y)-y}$ which does not depend on $y \in S$.
Let us call this bound $\mu/2$. Then:
\[
\nr{y-x}_1 - \mu \leq \nr{\pi(y) - \pi(x)}_1 \leq \nr{y-x}_1 + \mu
\]
And so, applying this inequality with fixed $x$ and all $y \in \tA(x,n)$, we
have:
\[
\bigg( \Bron(\pi(x),n-\mu) \cap \Pron(x) \bigg) \subseteq A(x,n) \subseteq
\bigg( \Bron(\pi(x),n+\mu) \cap \Pron(x) \bigg)
\]
where $\Bron$ stands here for the ball in $E$.

A consequence is that, given a patch $A(x,n)$, the number of different
$y \in \tA(x,n)$ such that $A(x,n) = A(y,n)$ is bounded, and its bound doesn't
depend on $x$.
Indeed, if $\nr{x-y} > 2 \mu$, then $\Bron(\pi(x),n-\mu) \cap \Pron(x)$ is
not included in $\Bron(\pi(y),n+\mu) \cap \Pron(y)$, so that we cannot have
$A(x,n) = A(y,n)$.
So a bound of the number of candidates for being the center of a given patch
is, for example, $\max_u{\big(\Pron(u) \cap \Bron(\pi(u),2\mu)\big)}=:M$,
which is finite by uniform discreteness of the Delone sets associated to the
cut and projection data.
This shows that a patch of size $n$ can have at most $M$ centers.
\ep\medbreak

\proc{Demonstration of the proposition.}
The proposition is an immediate consequence of the lemma above: as a given
patch $A(x,k)$ has at most $M$ centers, it means that such a patch corresponds
to at most $M$ pointed patches $A(x,k)$ with distinguished point
$y \in A(x,k)$ ($y$ being a center for the patch).
Therefore, we have the result with $\lambda = 1 / M$.
The other inequality is trivial.
\ep\medbreak

\subsection{The Structure of Singular Points}

In the next sections of this paper, the singular points will play a crucial role
to compute the complexity function. Therefore, it is necessary to understand
their geometrical structure.

\begin{prop}
\label{prop-struct-pts-k-sing}
We have the following properties:
\begin{equation}\label{eqn1}
\sg = \partial S + \Z^N
\end{equation}
and
\begin{equation}\label{eqn2}
\sg_k = \big( \Bron (0,k) \cap \Z^N + \partial S \big) \cap S
\end{equation}
so that $\sg_k$ is an increasing sequence of subsets of $S$, the union of
which is $\sg \cap S$.
Moreover, $\sg_k$ and $\sg$ are invariant under the action of $E$ by
translation.
\end{prop}

\proc{Demonstration.}
It is straightforward, using the definition.
\ep\medbreak

As the singular points are invariant by the action of $E$, it is enough to
study their restriction to $F$, which is a complementary subspace of $E$ in
$\R^N$.
In the next sections, we will mostly work on $F$.
Applying $\pi_F$ to the equations above, we can see that the first condition
which defines an almost-canonical acceptance domain in definition~\ref{hypK} is
equivalent to the fact that the singular points are a countable union of
hyperplanes.

We will state two results on singular points. The first one will be a
qualitative result, the second one will be a quantitative (and a bit
technical) result.
The proofs will be quite straightforward in the case of almost-canonical
acceptance domains, as the definition of almost-canonical was specifically
designed so that these results are true.
In order to prove that they hold for canonical acceptance domains as well, we
prove the following proposition.
\begin{prop}
A canonical acceptance domain is almost-canonical.
\end{prop}

\proc{Demonstration.}
Consider a canonical acceptance domain $K$. Let $f$ be a $(N-d-1)$-dimensional
face of $K$, and $H_f$ the hyperplane parallel to it. Up to renaming the
vectors, the face $f$ is of the following form:
\[
 f = v + \left\{ \sum_{k=1}^{N-d-1}{\lambda_k \pi_F(e_k)} \tq
0 \leq \lambda_k \leq 1 \right\}
\]
where $v \in V$.
Note $\Lambda$ the group generated by the $(\pi_F (e_i))_{i=1}^{N-d-1}$.
This group is a regular lattice of $H_f$ and a subgroup of $H_f \cap \Gamma$.
Furthermore, $f-v$ it is a fundamental domain for the action of $\Lambda$,
Therefore, $f + \Gamma$ is a union of hyperplanes. The same argument applies to
all the faces, and $\partial K + \Gamma$ is a union of hyperplanes.

To prove the second condition, remark once again that $f-x$ is compact
for any $x \in f$, and is a fundamental domain for the action of the locally
finite group $\Lambda$.
Therefore, there exist only a finite number of translates of the form
$f-x + \gamma$ with $\gamma \in \Lambda$ that intersect a given
neighborhood of the origin. As $f$ is bounded, this finite number can be chosen
such as not to depend on $x$.
On the other hand, $f-x + \Lambda$ covers $H_f$.
So a given neighborhood of the origin, it is covered by a finite number of
translates of $f-x$ of the form above.
It proves that a canonical acceptance domain is almost-canonical.
\ep\medbreak

Let us now make precise the structure of singular points.
\begin{definition}
We note $H_1, \ldots, H_m$ all the distinct hyperplanes of $F$ which
are parallel to faces of $K$.
We call these hyperplanes the \emph{singular hyperplanes}, or \emph{cut
hyperplanes} of $F$.
Remember that $\Gamma := \pi_F (\Z^N)$. We define
$\Gamma^i := \Stab_\Gamma (H_i)$.
\end{definition}

Remark that as a face $f$ parallel to $H_i$ is compact, and $f + \Gamma^i$ is
an hyperplane, then $\Gamma^i$ always contains a regular lattice of $H_i$. It
implies that it is at least of rank $N-d-1$.
So the hyperplanes $H_i$ are always directed by $N-d-1$ vectors of $\Gamma$.

\begin{figure}[htp]
\begin{center}
\includegraphics[scale=0.6]{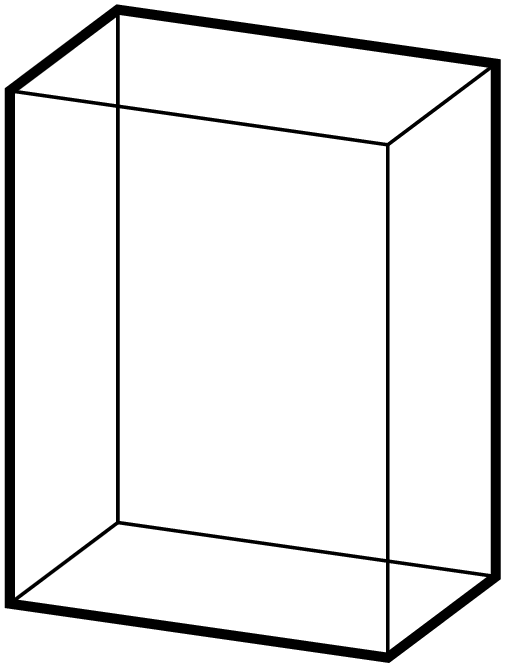} \qquad \quad
\includegraphics[scale=0.6]{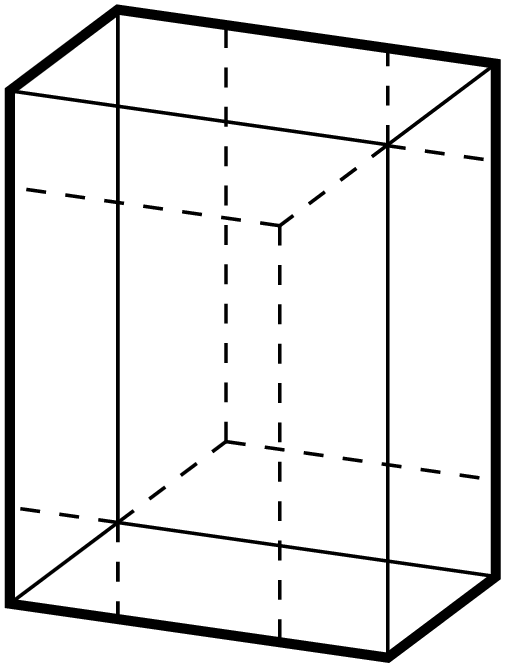} \qquad \quad
\includegraphics[scale=0.6]{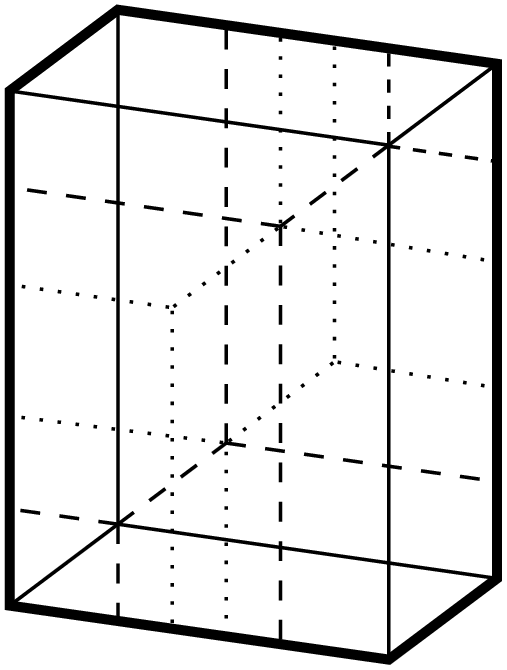}
\caption{Respectively the $1$, $2$ and $3$-singular points
of a $(3,1)$-cut and projection tiling, seen on $F$ (for a canonical acceptance
domain)}
\label{fig-sg-pts}
\end{center}
\end{figure}

The following proposition follows quite easily from definition~\ref{hypK}.
It is best understood with a picture: see figure~\ref{fig-sg-pts}.
\begin{prop}
\label{prop-hyp-sing}
The singular points of $F$ are the union of translates of singular hyperplanes:
\[
\sg \cap F = \bigcup_{i=1}^m{\left(H_i + V_i \right)} + \Gamma
\]
where $V_i$ is the subset of $V$ consisting of the vertices of $K$ which belong
to a face parallel to $H_i$.
In particular, we always have:
\[
\sg \cap F \subseteq \left(\bigcup_{i=1}^m{H_i} \right) + V + \Gamma
\]
and if $V \subseteq \Gamma$ (for example in the case of a canonical acceptance
domain), we have:
\[
\sg \cap F = \left( \bigcup_{i=1}^m{H_i} \right) + \Gamma
\]
\end{prop}

The following technical lemma is a qualitative version of the previous
proposition. It will use the second property of almost-canonical acceptance
domains (cf. definition~\ref{hypK}).
Once again, having in mind figure~\ref{fig-sg-pts} helps.
\begin{lemma}
\label{lemma-hyp-sing}
We use the notations above.
Then for all $1 \leq i \leq m$, there exists $\eps_i > 0$ and $M_i \in \N$ such
that for all $y \in (\pi_F(\Bron(0,n) \cap \Z^N) + V_i)$, the points of
$(H_i + y) \cap \Bron(y,\eps_i)$ are $(n+M)$-singular.
\end{lemma}

\proc{Demonstration.}
First, remark that if we apply $\pi_F$ to equation~(\ref{eqn2}), we have:
\[
\sg_n \cap F = \partial K + \pi_F \big( \Z^N \cap \Bron(0,n) \big)
\]
Let $1 \leq i \leq m$.
Consider $\eps_i > 0$ such that, according to definition~\ref{hypK}, the ball of
radius $\eps_i$ around $0$ in $H_i$ is covered by a finite number of translates
of any face parallel to $H_i$ (there are only finitely many such faces).
Let $n \in \N$ and $y \in \pi_F(\Bron(0,n) \cap \Z^N) + V_i$.
We write $y = v + \pi_F (z)$ with $z \in \Z^N$ such that $\nr{z}_1 \leq n$, and
$v \in V_i$.
Call $f$ the face parallel to $H_i$ to which belongs $v$.
But there are a finite number of translates of $f - v$ by elements of
$\Gamma^i$ which cover $\Bron (0, \eps_i)$.
Let us call $\gamma_1, \ldots \gamma_k$ these elements, and let $M_i$ be:
\[
 M_i = \max_{1 \leq i \leq k}{\min\{\nr{w}_1 \tq w \in \Z^N,
\pi_F (w) = \gamma_i\}}
\]
Note that $M_i$ can be chosen such as not to depend on $v$ (as $f$ has only
finitely many vertices).
So if $x \in \Bron(y,\eps_i) \cap H_i$, then
$x-y \in \Bron (0,\eps_i) \cap H_i$, and $x-y$ is in some
$f - v + \pi_F (w)$ with $w \in \Z^N \cap \Bron (0,M_i)$.
It means that $x \in f + \pi_F(z+w)$ with $z+w \in \Z^N \cap \Bron(0,n+M)$.
\ep\medbreak

We know (proposition~\ref{prop-hyp-sing}) that if a point $y \in \Gamma + V_i$
is $n$-singular, then $H_i + y$ is singular.
What this lemma states is that, restricted to $H_i$, there exists a
neighborhood of $y$ for the relative topology of $H_i$ which is
$(n+M)$-singular, and $M$ does not depend on $n$.
It states that locally, around a $n$-singular point, the $(n+M)$-singular
points look like hyperplanes.

Therefore, we have two inclusions: the set of $n$-singular points is included in
a union of hyperplanes, and it contains (at least locally) a union of
hyperplanes.

\section{Complexity for Cut and Projection Tilings}

\subsection{Statement of the Theorem}

\begin{figure}[htp]
\begin{center}
\includegraphics[scale=0.70]{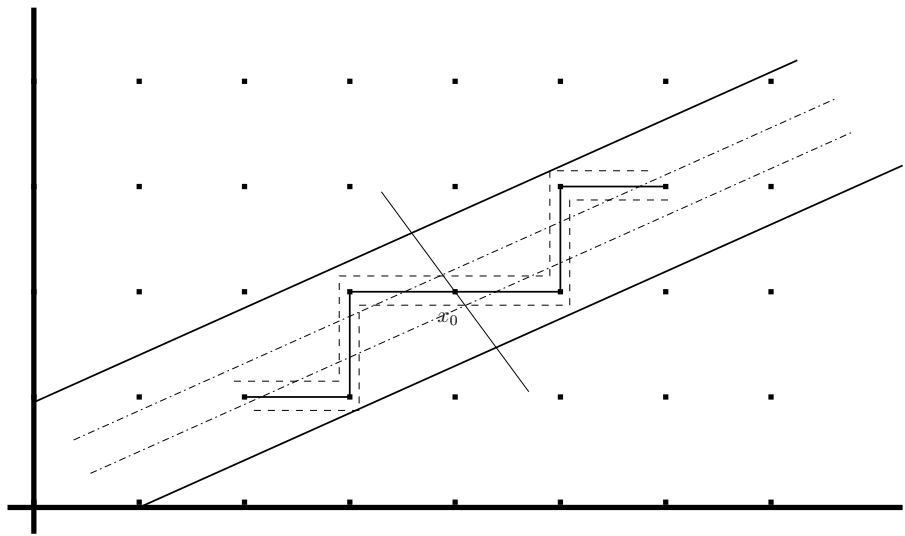}
\label{fig-cc}
\caption{When we move $x_0$, the seed of the lifted patch, we do not
change its shape, until it becomes $3$-singular (the patch has size $3$).
As long as $x_0$ stays in the band (the connected component of $\reg_3$,
here pictured with dash-dotted lines), the patch $\tA (x_0,3)$ will not change
up to translation.}
\end{center}
\end{figure}

We will state in this section the key proposition and the main theorem of this
paper.
The first proposition states that, up to a constant, to count the number
of different patches of given size $n$ amounts to count the number of
connected components of $\reg_n \cap K$. It is the precise statement
corresponding to the fact that when considering the patch $A(x,n)$,
$\pi(x)$ encodes its position in $E$, while $\pi_F(x)$ encodes its shape.
This is illustrated by figure~\ref{fig-cc}.
\begin{prop}
\label{prop-fondam}
We note $p_\mathrm{pt}$ the pointed complexity function, as defined above.
For all $n \in \N \cup \{0\}$, $c(n)$ will denote the number of connected
components of $\reg_n$.
Then, for all $n \in \N$, we have the following.
\[
p_\mathrm{pt}(n) = c(n)
\]
\end{prop}

The following corollary results immediately from the proposition above and
proposition~\ref{prop-center}.
\begin{cor}
We note $p$ the complexity function and $c$ the function which counts the
number of connected components of $\reg_n$.
Then there exists $0 < \lambda < 1$ such that for all $n \in \N$, we have:
\[
\lambda c(n) \leq p(n) \leq c(n)
\]
\end{cor}

We will now state the main theorem of this paper, which makes a link between
the complexity of cut and projection tilings associated to given data, and
the rank of the stabilizer in $\Gamma$ of the singular hyperplanes.

\begin{theo}
\label{theorem}
Let $(\R^N = E \oplus F, \Z^N, \pi, K)$ be the data of a cut and projection
method with $K$ almost-canonical.
Let $p$ be the associated complexity function.
Let $(H_i)_{i=1}^m$ be the singular hyperplanes in $F$, and
$\Gamma := \pi_F (\Z^N)$.
We define the following objects:
\begin{itemize}
\item for all $1 \leq i \leq m$, let $\Gamma^i := \Stab_\Gamma (H_i)$;
\item let $\alpha_i + 1$ be the rank of $\Gamma / \Gamma^i$.
\item let $(\alpha_i)_{i \in I}$ be $N-d$ elements amongst the
$(\alpha_i)_{i=1}^m$ satisfying $\bigcap_{i \in I}{H_i} = \{0\}$, and such
that $\alpha := \sum_{i \in I}{\alpha_i}$ is maximal.
\end{itemize}
Then there exist $C_1, C_2 > 0$ such that for all $n \in \N$,
\[
C_1 n^\alpha \leq p(n) \leq C_2 n^\alpha
\]
\end{theo}

\begin{rem}
Generically, the rank of $\Gamma$ is $N$, and the rank of the stabilizers of
the $H_i$ is minimal, which means it is equal to $N-d-1$. So for all $i$,
$\alpha_i = d$, and the complexity $p(n)$ grows like $n^{d(N-d)}$.
This is also the maximal complexity of a $(N,d)$-cut and projection tiling
with our assumptions on $K$ (in particular, with canonical acceptance domain).
It is easily seen that $alpha \geq N-d$. There is another less trivial lower
bound for $\alpha$:
in theorem~\ref{theorem-cohom}, we will prove that $\alpha = d$ is the
minimal value of $\alpha$ for a cut and projection Delone set with no period.
\end{rem}

Proposition~\ref{prop-fondam} shows that the functions $p(n)$ (complexity)
and $c(n)$ (number of connected components of $K \cap \reg$) have the same
order of growth.
The theorem is then proved by counting the number of connected components of
$K \cap \reg$, using a few geometric lemmas, which we will state and prove
in the following sections.

\proc{Demonstration of proposition~\ref{prop-fondam}}
Let us first show that for all $n \in \N$, $p_\mathrm{pt}(n) \leq c(n)$.
Remember that $p_\mathrm{pt}(n)$ is the number of pointed patches
$A_{\mathrm{pt}}(x,n)$ which differ up to a translation by an element of
$E$, with $x \in K \cap \reg_n$.
What we have to show is that, if $x$ and $y$ are in the same connected
component of $\reg_n \cap K$, then $(A(x,n),\pi(x)) = (A(y,n),\pi(y))$.
Let $x \in \reg_n \cap K$. We will show that the set $X$ of all the
$y \in \reg_n \cap K$ such that $A(x,n) = A(y,n)$ and $\pi(x) = \pi(y)$ is both
open and closed for the relative topology on $\reg_n \cap K$.
First, remark that it is superfluous to check $\pi(x) = \pi(y)$, as by
hypothesis, $x,y \in K$, and so $\pi (x) = \pi (y) = 0$.
Let $y \in X$, and $d$ be the distance between $\tA (y,n)$ and $\partial S$,
so that $\Bron (0,d/2) + \tA (y,n)$ is included in the interior of $S$.
Then $\Bron (y,d/2) \cap K \subseteq X$, and $X$ is open.
Conversely, consider $(y_k)_{k \geq 0}$, a sequence in $X$ which converges to
$y \in \reg_n \cap K$.
Let $d$ be the distance between $\tA(y,n)$ and $\partial S$.
Then for all $y' \in \Bron(y,d/2)$, $A(y,n) = A(y',n)$. But for $k$ big enough,
$y_k \in \Bron(y,d/2)$, so that $A(y,n) = A(y_k,n) = A(x,n)$, and $y \in X$.
The set $X$ being open and closed, it contains the connected component of
$\reg_n$ containing $x$, and $p_\mathrm{pt}(n) \leq c(n)$.

Let us now show the other inequality.
Let $A_\mathrm{pt} (x,n)$ be a pointed patch of size $n$; we will show
that there exists exactly one connected component $\Cron$ such that if
$y \in C$, $A_\mathrm{pt}(y,n) = A_\mathrm{pt}(x,n)$ up to translation.
If $A_\mathrm{pt}(x,n) = A_\mathrm{pt}(y,n)$, by injectivity of
$\restr{\pi}{\Z^N}$, we have that $\tA(x,n)$ and $\tA(y,n)$ differ by an
element of $F$.
Moreover, as we consider pointed patches, this element is exactly $y-x =: v$.
Then we have $\tA(x,n) + v = \tA (y,n)$.
Both of these sets being included in $S$, by convexity, we have that
$\tA(x,n) + tv = \tA(x+tv,n)$ is included in $S$ for all $t \in [0,1]$, so that
all the $x+tv$ are $n$-regular, and $x$ is in the same connected component of
$\reg_n$ as $x+v = y$.
Thus, $x$ and $y$ are in the same connected component of $\reg_n$.
\ep\medbreak

\begin{rem}
The second part of the demonstration above proves that the connected components
of $\reg_n$ are convex.
\end{rem}

We will leave the proof of theorem~\ref{theorem} for the next sections, but let
us give a motivation for the technical lemmas we are going to state (and
prove).
Our goal is to give an asymptotic estimation of $c(n)$, the number of
connected components of $\reg_n \cap K$. Remember that the set of regular
points is a convex polytope $K$ minus hyperplanes:
\[
\reg \cap K = K \setminus \left( \Gamma + \bigcup_{i=1}^m{(H_i + V_i)} \right)
\]
where the $V_i$ are finite sets.
Similarly, one can prove that:
\[
 \sg_n \subseteq
\Big( \big( \pi_F (\Z^N \cap \Bron(0,n)) + V \big) \cap K \Big) +
\bigcup_{i=1}^m{H_i}
\]
The idea is to compare $c(n)$ to the number of connected
components of $K$ minus a finite number of hyperplanes in each direction.
To do so, we have to count the number of hyperplanes actually cutting $K$ in the
formula above.
For this, according to the formula above, we might need to count the number of
points of $\pi_F (\Z^N \cap \Bron(0,n)) + V$ which fall inside of $K$.
Actually, the set we will count will not be exactly this one, as we will have
to take into account the fact that some hyperplanes are stabilized by elements
of $\Gamma$. However, the idea remains the same.
Then, we will need to estimate the number of connected components of a compact
convex subset of $\R^{N-d}$ when it is cut by hyperplanes.
We will state two lemmas which will answer these questions.
Combining these two lemmas, we will then prove the theorem.

\subsection{Preliminary Lemmas}

We begin by an ``equidistribution-like'' lemma: this lemma aims at giving an
estimation of the number of points of $\pi_F(\Bron(0,n) \cap \Z^N) \cap K$,
when $n$ tends to infinity.

\begin{lemma}\label{lemma1}
Let $F$ be a $p$-dimensional vector space, $\gamma_1, \ldots, \gamma_q$ be
$q$ vectors of $F$, which span a dense subgroup of rank $q$ ($q > p$).

Then for all bounded sets $U$ with non empty interior, there exist $C_1,C_2>0$,
such that for all $n$, we have:
\[
C_1 n^{q-p} \leq \Card{\bigg\{m = (m_i)_{i=1}^q \in \{1, \ldots, n\}^q \tq
\sum_{i=1}^q{m_i \gamma_i} \in U \bigg\}} \leq C_2 n^{q-p}
\]
\end{lemma}

We can apply this lemma with $p = N-d$, $q = \rk{\Gamma}$, and where
the $\gamma_i$ are a system of generators of $\Gamma$ chosen
amongst the $(\pi_F (e_i))_{1 \leq i \leq N}$.
In the special case where we make the additional assumption that $\Z^N \cap
E = \{ 0 \}$, we will have $q = N$. In this case, this lemma states
that the cardinality of $\pi_F \big( \tA (0,n) \big)$ grows like $n^d$.
This gives the cardinality of a patch of size $n$.
An interpretation is that the number of integer points contained in the
convex set $(K+E) \cap B(0,n)$ grows like the volume of this set, which
is roughly speaking the volume of $K + (E \cap B(0,n))$, that is, grows like
$n^d$, up to constants.

\proc{Demonstration.}
Let us first remark that by precompacity of bounded sets, it is enough to check
the property for a set $U$ of the form $U = \Bron(0, \eps')$ with small $\eps'$.

Up to re-ordering the vectors, we assume that $\gamma_1, \ldots, \gamma_p$
define a parallelotope of maximal volume amongst all the parallelotopes which
can be formed from the $\{\gamma_i\}_{i=1}^q$.
Thus, $(\gamma_i)_{i=1}^p$ form a basis of $F$.
Let $\nr{.}_2$ and $\nr{.}_\infty$ be respectively the Euclidean and sup norms
associated to this basis.
One can check that for all $i > p$, we have $\nr{\gamma_i}_\infty \leq 1$.
Let $\Lambda$ be the lattice of $F$ generated by $\gamma_1, \ldots, \gamma_p$.

We will prove a result of equidistribution relative to the subgroup of the
torus $T:= F / \Lambda$ generated by the vectors $(\gamma_i)_{i>p}$.
Let $i > p$. By hypothesis, the subgroup of $E$ generated by all the
$(\gamma_j)_{j=1}^q$ is of rank $q$, so that the group $\Z \gamma_i$ is
infinite as a subgroup of $T$.
Let us call $T_i$ its adherence. This is a sub-torus of $T$, in which
$\Z \gamma_i$ is dense.
It is then classical (see for example~\cite{Kor}), that the sequences $(n
\gamma_i)_{n\geq 0}$ and $(-n \gamma_i)_{n>0}$ are equidistributed in $T_i$.
Thus, for all $\eps > 0$, if $V(\eps)$ is a small neighborhood of $0$ in $T$
(which is the image on the quotient of a small ball $\Bron(0,\eps)$ in $F$,
with $\Bron(0,\eps)$ being included in a fundamental domain),
then $V(\eps) \cap T_i$ has positive measure (for the invariant measure on
$T_i$), and thus:
\[
\frac{1}{2n+1} \Card \Big\{
  k \in \{-n \ldots n \} \tq k \gamma_i \in V(\eps) \subseteq T \Big\}
  \xrightarrow[{n \rightarrow +\infty}]{} l_i > 0
\]
Lifting this result, we deduce that:
\begin{multline*}
\frac{1}{2n+1} \Card \Big\{
  k=(k_j)_{j=1}^{p+1} \in \Z^p \times \{-n \ldots n \} \tq \\
  \sum_{j=1}^{p}{k_j \gamma_j} + k_{p+1} \gamma_i \in \Bron(0,\eps) \subseteq E
\Big\}
  \xrightarrow[{n \rightarrow +\infty}]{} l_i > 0
\end{multline*}
To simplify the notation, let us define for all $i>p$:
\[
A^{(i)}_n(\eps) := \Big\{ k=(k_j)_{j=1}^{p+1} \in \{-n \ldots n \}^{p+1} \tq
  \sum_{j=1}^{p}{k_j \gamma_j} + k_{p+1} \gamma_i \in \Bron(0,\eps) \subseteq E
\Big\}
\]
Remark that, the condition on the norms (that is $\nr{\gamma_i}_\infty \leq 1$)
implies that $\nr{n \gamma_i}_\infty \leq \abs{n}$ for all $n \in \Z$.
Thus, for all $n \in \Z$, there exists a unique vector $v$ such that
$n \gamma_i + v$ belongs to a fundamental domain containing $\Bron (0,\eps)$ of
the lattice $\Lambda$, and we have $\nr{v}_\infty \leq n$.
Thus, we deduce the following:
\[
\frac{1}{2n+1} \Card A^{(i)}_n(\eps) \xrightarrow[{n \rightarrow +\infty}]{}
l_i > 0
\]
Remark that the same reasoning holds for all $\gamma_i$, with $i > p$.

Let us now prove the lemma itself: let $\eps' > 0$, we want to find both an
upper bound and a positive lower bound for the sequence
$(2n+1)^{-(q-p)} \Card{A_n(\eps')}$, where
\[
A_n(\eps') := \big\{k=(k_j) \in \{-n, \ldots, n \}^q \tq
   \sum_{j=1}^{q}{k_j \gamma_j} \in \Bron(0,\eps') \subseteq E \big\}
\]
Let us first find bounds to the cardinality of $(2n+1)^{-(q-p)} A'_n(\eps')$,
where
\begin{multline*}
A'_n(\eps'):=\Bigg\{ k=(k_j)_{j=1}^q \in \big\{-(q-p)n, \ldots, (q-p)n \big\}^p
   \times \{-n, \ldots, n \}^{q-p} \tq   \\
   \sum_{j=1}^{q}{k_j \gamma_j} \in \Bron(0,\eps') \subseteq E \Bigg\}
\end{multline*}
Consider $k^{(i)} \in A^{(i)}_n (\eps)$, for all $i \in \{p+1, \ldots, q\}$,
with $\eps := \eps'/(q-p)$. Then to these $(k^{(i)})_{i=p+1}^q$, one can
associate a $k \in A'_n(\eps')$ defined by $k_j := \sum_{i=p+1}^q{k^{(i)}_j}$ if
$j \leq p$, and $k_j := k^{(j)}_{p+1}$ if $j > p$.
Plus, this application
\[
\prod_{i=p+1}^q{A^{(i)}_n(\eps)} \longrightarrow A'_n(\eps')
\]
is one-to-one.
Indeed, if $(k^{(i)})_{i=p+1}^q \neq (k'^{(i)})_{i=p+1}^q$ (with respective
images $k$ and $k'$ in $A'_n (\eps')$), then there exists
$p+1 \leq j \leq q$ such that $k^{(j)} \neq k'^{(j)}$.
It necessarily means that $k^{(j)}_{p+1} \neq k'^{(j)}_{p+1}$, because
if  $k^{(j)}_{p+1} = k'^{(j)}_{p+1}$, then as $\eps'$ (and thus $\eps$) is 
chosen small enough, there exists a \emph{unique} vector
$v \in \{-n, \ldots, n\}^p$ such that
$k^{(j)}_{p+1} \gamma_j + v \in \Bron (0,\eps)$.
Thus, $k^{(j)} = (v,k^{(j)}_{p+1}) = k'^{(j)}$, which is a contradiction.
Thus $k^{(j)}_{p+1} \neq k'^{(j)}_{p+1}$, and so $k \neq k'$, and the
application defined above is one-to-one.

As a consequence, the cardinality of $A'_n(\eps')$ is bigger than the product of
the cardinalities of the $A^{(i)}_n(\eps)$, and so:
\[
\liminf_{n\rightarrow \infty}{\frac{1}{(2n+1)^{q-p}} \Card A'_n(\eps')} > 0
\]
Let us now find an upper bound to the cardinality of the $A'_n(\eps')$: if
$\eps'$ is chosen small enough, \ie $\Bron (0, \eps')$ is included in a
fundamental domain of the lattice $\Lambda$, then:
\[
\Card (A'_n(\eps')) \leq (2n+1)^{q-p}
\]
because an element of $A'_n(\eps')$ is entirely determined by its $q-p$ last
coordinates.

Now, the cardinality of $A_n$ can be easily compared to the cardinality of
$A'_n$, as we have:
\[
\Card A'_n(\eps') \leq \Card A_{(q-p)n}(\eps') \leq \Card A'_{(q-p)n}(\eps')
\]
As the upper bound we want to find is polynomial in $n$, and
$\Card (A_n(\eps'))$ is increasing, we have the result.
\ep\medbreak

The next lemma gives an upper bound for the number of connected
components of a polytope cut by affine hyperplanes, assuming we have an
estimation of the number of hyperplanes involved.
As we said, $\reg_k$ is ``almost'' $K$ minus translates of the singular
hyperplanes. Therefore, this lemma will be a key lemma in order to have
an upper bound of $c(n)$.

\begin{lemma}\label{lemma2}
Let $F$ be a $p$-dimensional real vector space.
Then for all $m \in \N$, for all $H_1, \ldots, H_m$ distinct linear
hyperplanes of $F$, for all $(\alpha_i)_{i=1}^m$ integers, for all finite
families of integer-valued increasing sequences $\big(\beta^{(i)}\big)_{i=1}^m$
such that for all $i$, $\beta^{(i)}_n = O(n^{\alpha_i})$, we have the
following:

If $K$ is an open convex and bounded subset of $F$ and $(K_n)_{n \geq 0}$ is
a decreasing sequence of sets, such that $K_0 = K$ and $K_{n}$ can be obtained
from $K$ by removing $\beta^{(i)}_n$ distinct translates of the $H_i$ (for all
$i$), then the number $c(n)$ of connected components of $K_n$ satisfies:
\[
c(n) = O (n^\alpha)
\]
where $\alpha := \sum_{i \in I}{\alpha_i}$, with $I$ defined as follows:
\begin{itemize}
\item $\bigcap_{i \in I}{H_i}$ is of minimal dimension;
\item amongst such sets, $I$ is of minimal cardinality (the cardinality of $I$
is then $p - \dim{\bigcap_{i=1}^m{H_i}}$);
\item amongst such sets, $I$ is such that $\alpha$ is maximal.
\end{itemize}
\end{lemma}

\proc{Demonstration.}
We prove this lemma by recurrence on $p$. If $p = 1$, then $m = 1$ as the only
hyperplane of a $1$-dimensional space is $\{0\}$. Then $K_n$ is a
segment with $\beta_1$ points removed, \ie a union of $\beta_1 + 1$ intervals.
So its number of connected components is $\beta_1 + 1 = O (n^{\alpha_1})$, and
the lemma is proved for $p=1$.

Let us assume that the lemma is true for all $(p-1)$-dimensional vector spaces.
Let $F$ be a $p$-dimensional vector space. Let $K$, $m$, $(H_i)_{i=1}^m$ 
$\big(\beta^{(i)}\big)_{i=1}^m$ and $(K_n)_{n \geq 0}$ be as in the statement
of the lemma.
We note $\Cron$ the function which associates its number of connected
components to a subset of $F$ ($\Cron$ takes values in $\N \cup \{\infty\}$),
so that $c(n) = \Cron (K_n)$.
To fix the notations, we note for all $i$, $\left(X^{(i)}_n\right)_{n \geq 0}$
an increasing sequence of finite subsets of $F$, so that
$\Card{X^{(i)}_n} = \beta^{(i)}_n$ for all $i,n$, and
\[
K_n = K \setminus \Big( \bigcup_{i=1}^m{(X^{(i)}_n + H_i)} \Big)
\]
We note
$c(x,k,n):=\Cron \bigg( (H_k + x) \cap \Big( K \setminus \bigcup_{i=1}^{k-1}
{(X^{(i)}_n + H_i)} \Big) \bigg)$,
and we remark that the connected components of $K_n$ are convex (being the
intersection of half-spaces), so that $c(x,k,n)$ is exactly the number of
components intersected by the hyperplane $H_k + x$.
Therefore, for all $x \in X^{(k)}_n$, we have:
\begin{equation}\label{eqn}
\Cron \bigg( K \setminus \Big( (x+H_k) \cup
\bigcup_{i=1}^{k-1}{(X^{(i)}_n + H_i)} \Big) \bigg)
=
\Cron\bigg(K \setminus \bigcup_{i=1}^{k-1}{(X^{(i)}_n + H_i)} \bigg)
+ c(x,k,n)
\end{equation}
An interpretation of the equation above is the following: if we have an open
set $K'$, made of a certain number $\Cron (K')$ of convex connected components,
we want to count the number of connected components of $K'$ minus an affine
hyperplane $H$.
The number of additional components obtained by cutting by this hyperplane
equals the number of components of $K'$ which are intersected by $H$, that
is the number of components of $K' \cap H$ (by convexity of the components
of $K'$).
Equation~(\ref{eqn}) expresses this fact with a particular $K'$, which allows
to iterate the process: if we consider that $K_n$ is nothing more than
$K$ with a certain number of hyperplanes removed \emph{one after another},
we can deduce $c(n)$ by successive sums and we have:
\[
c(n) = \sum_{k=0}^m{\Bigg( \sum_{x \in X^{(i)}_n}{c(x,k,n)} \Bigg)}
\]

So if we find an upper bound of the number of $c(x,k,n)$ independently of
$x$ for all $k$, we should be able to find an upper bound of $c(n)$.
The upper bound of $c(x,k,n)$ will be obtained by using the recurrence
hypothesis.

Let $k \in \{1, \ldots, m\}$. We want to have an upper bound of $c(x,k,n)$.
Consider the following set:
\[
\left( K \setminus \bigcup_{i=1}^{k-1}{(X^{(i)}_n + H_i)} \right) \cap (H_k+x)
\]
This set is obtained by considering the space $F':=H_k + x$ (which we regard as
a vector space of dimension $p-1$ by taking $\{x\}$ as the origin), and a
compact convex subset, $K':= K \cap F'$ which is cut by hyperplanes directed by
the $H_j \cap H_k$, where $j \neq k$.
It can happen that some of these hyperplanes are not distinct, so we rename
them $H'_1, \ldots, H'_s$, with $s < m$.
The sequence of subsets $(K_n)_{n \geq 0}$ of $K$ induces a sequence
$(K'_n)_{n \geq 0}$ of subsets of $K'$, obtained by removing hyperplanes
directed by the $H'_i$.
The number of hyperplanes removed at each step is given by sequences
$(\gamma^{(i)}_n)_{n \geq 0}$ for $1 \leq i \leq s$.
Let us now find an estimation of the $(\gamma^{(i)}_n)_{n \geq 0}$ to apply the
recurrence hypothesis: if $H'_i = H_{i(1)} \cap H_k = \ldots =
H_{i(q)} \cap H_k$, then $\gamma^{(i)}_n = O \big(n^{\delta_i} \big)$, where
$\delta_i := \max_j{(\alpha_{i(j)})}$.
The fact that we required the $(H_i)_{i \in I}$ to be transverse (\ie
$\bigcap_{i \in I}{H_i}$ of minimal dimension) implies that if $i,j \in I$,
then $H_k \cap H_i \neq H_k \cap H_j$.
Therefore, the sum of the $p-1$ biggest elements amongst the
$\{\delta_i\}_{i=1}^s$ (we call this sum $\delta$) with suitable additional
condition of transversality for the associated hyperplanes, is lesser than
the sum of the $p-1$ biggest elements amongst the $\{\alpha_i\}_{i=1}^m$ (with
the same condition of transversality).

So by recurrence, we have $c(x,k,n) = O(n^\delta)$.
This upper bound does not depend on $x$ so if we sum the $x \in X^{(k)}_n$, we
get that $\sum_{x \in X^{(i)}_n}{c(x,k,n)} = O (n^{\delta + \alpha_k})$.
But $\delta + \alpha_k \leq \alpha$. We deduce that
$\sum_{x \in X^{(i)}_n}{c(x,k,n)} = O (n^{\alpha})$ and so, as we can have the
same bound for all $k$, summing again, we obtain the expected result:
\[
c(n) = O (n^\alpha)
\]
\ep\medbreak

\subsection{Proof of the Theorem}

Let us now prove the theorem~\ref{theorem}. The proof falls in two parts:
upper, and lower bound.
Actually, thanks to proposition~\ref{prop-fondam}, the quantity we will bound
(above and below) is $c(n)$ rather than $p(n)$.
We will then use lemmas~\ref{lemma-hyp-sing}, \ref{lemma1} and \ref{lemma2}
to count the number of connected components.

\proc{Upper bound.}
Remember that $V$ is the set of vertices of $K$, seen as a polygon.
As $\Gamma$ is a free Abelian group of rank $r$, we identify $\Gamma$ and
$\Z^r$ by the choice of a group-basis $\gamma_1, \ldots, \gamma_r$, such that
the $\gamma_i$ are projections of elements of the canonical basis of $\R^N$.
If $a \in \Gamma$, and $a = (a_i)_{i=1}^r$ when decomposed in this basis, we
define:
\[
\Bron_\Gamma (a,n) := \{ b = (b_i)_{i=1}^r \in \Gamma \tq
\max_i{\abs{b_i-a_i}} \leq n \}
\]
There exists a constant $\rho$ such that
$\pi_F (\Z^N \cap \Bron(0,1)) \subseteq \Bron_\Gamma (0,\rho)$.

We already mentioned that:
\[
\sg_n \cap K = \Big(\pi_F \big(\Bron(0,n) \cap \Z^N \big) + \partial K \Big)
\cap K
\]
(apply $\pi_F$ to equation~(\ref{eqn2})).
As $\partial K$ is included in $V + \bigcup_{i=1}^m{H_i}$, we have:
\[
\sg_n \cap K \subseteq
\bigg( \Bron_\Gamma(0,\rho n) + V + \bigcup_{i=1}^m{H_i} \bigg) \cap K
\]
Given $i \in \{1, \ldots, m\}$, we now have to find an upper bound of the number
of distinct hyperplanes of the following form:
\[
\Bron_\Gamma(0,\rho n) + V + H_i
\]
and which actually intersect $K$. We call this number $\beta^{(i)}_n$.
We have to determine the asymptotic behavior of the sequence
$(\beta^{(i)}_n)_{n \in \N}$.
We will then apply lemma~\ref{lemma2}.

Let $q_i : F \rightarrow F / H_i$ be the quotient map. It maps
$F$ to a one dimensional $\R$-vector space, and it maps $\Gamma$ to
$\Gamma / \Gamma^i$, which is of rank $\alpha_i + 1$.
We want to find an upper bound of the number of hyperplanes of the form
$H_i + x$, with $x \in V + \Bron_\Gamma(0,\rho n)$ intersecting $K$.
We write these hyperplanes $H_i + \gamma + v$ with
$\gamma \in \Bron_\Gamma(0,\rho n)$ and $v \in V$.
But for a given $v \in V$, it amounts to count the number of points of
$q_i\big( \Bron_\Gamma(0,\rho n) + v \big)$ which belong to $q_i (K)$.
Equivalently, it amounts to count the number of points of 
$q_i\big( \Bron_\Gamma(0,\rho n) \big)$ which belong to $q_i (K - v)$.
Choosing a basis for $q_i(\Gamma) = \Gamma / \Gamma^i \simeq \Z^{\alpha_i + 1}$,
there exists a $\rho' > 0$ such that $q_i (\Bron_\Gamma(0,\rho n))
\subseteq \Bron_{\Gamma / \Gamma^i}(\rho' n)$ at least for all $n$ big enough.
We now apply lemma~\ref{lemma1}: the number of points in the set
$\Bron_{\Gamma / \Gamma^i}(\rho' n)$ which fall in $q_i (K - v)$ is dominated
by $n^{\alpha_i}$.
We do the same reasoning for for all $v \in V$, we sum the inequalities,
and we have:
\[
\beta^{(i)}_n = O (n^{\alpha_i})
\]

This reasoning holds for all $1 \leq i \leq m$, so $\big(\beta_n^{(i)}\big)_{n
\in \N}$ is dominated by $n^{\alpha_i}$ for all $i$.
Plus, as we said, $\beta^{(i)}_n$ is the number of hyperplanes parallel to $H_i$
of the form above intersecting $K$.
Then, we apply lemma~\ref{lemma2}, and we have:
\[
p(n) \leq c(n) = O (n^\alpha)
\]
\ep\medbreak

\proc{Lower bound.}
Let us then consider the set of indices $I$ satisfying condition that
$\alpha = \sum_{i\in I}{\alpha_{i}}$, and $\bigcap_{i \in I}{H_i} = \{0 \}$.
Let $\eps$ be the minimum over $i \in I$ of the $\eps_i$.
Consider an open parallelotope $P$ of non empty interior, of diameter
lesser than $\eps$, included in $K$, and the faces of which are parallel to the
hyperplanes $(H_i)_{i \in I}$.
Such a parallelotope does exist because all the $(H_i)_{i \in I}$ are
transverse.

Using the same reasoning as in the upper bound part above, we can say that, for
all $i \in I$ and for a given $v_i \in V_i$, the number
of points in $\big(\pi_F(\Bron (0,n \cap \Z^N)) + v_i\big) \cap P$ modulo
$H_i$ is greater than $C_i n^{\alpha_i}$ for a certain $C_i > 0$.

It implies, by lemma~\ref{lemma-hyp-sing}, that the $(n + M)$-singular points
inside $P$ contains at least $C_i n^{\alpha_i}$ distinct affine
hyperplanes (intersected with $P$) directed by $H_i$ for all $i \in I$.
Thus the number of connected components of the $(n+M)$-singular points, when
restricted to $P$, is at least $\prod_{i}{C_i n^{\alpha_i}}$, that is
at least $C n^\alpha$, for suitable $C$.
So we have:
\[
p(n) \geq \lambda c(n) \geq \lambda C n^\alpha
\]
\ep\medbreak

\section{Examples and Applications}

We now apply our results to three examples: the case $d=1$, the case $d=N-1$,
and the case of the octagonal tiling (a specific case with $d=2$ and $N=4$).

\subsection{The Octagonal Tiling}
The octagonal tiling can be defined as a cut and projection tiling with
$N = 4$ and $d = 2$. We denote $(e_1, \ldots, e_4)$ the canonical basis of
$\R^4$.

We consider the following matrix:
\[
A = \left( \begin{array}{cccc}
  0 & 0 & 0 & -1 \\
  1 & 0 & 0 & 0  \\
  0 & 1 & 0 & 0  \\
  0 & 0 & 1 & 0
\end{array} \right)
\]
Its characteristic polynomial is $X^4 + 1$, so that its invariant subspaces
are two orthogonal subspaces of $\R^4$ of dimension $2$.
Let us call $E$ and $F=E^\perp$ these subspaces ($E \perp F$ because the
matrix $A$ is orthogonal), and $\pi$ the orthogonal projection of
$\R^4$ on $E$.

We compute that we can choose the following basis for $F$:
\[
 F = \vect{\Big((\sqrt{2}/2,1/2,0,-1/2),(0, 1/2, \sqrt{2}/2, 1/2) \Big)}
\]
and that these vectors form an orthonormal basis of $F$. Thus, we can write the
matrix of $\pi_F$ in the canonical basis:
\[
\mat_{\mathrm{can}}{\pi_F} =
\left( \begin{array}{cccc}
  1/2         & \sqrt{2}/4 & 0          & -\sqrt{2}/4 \\
  \sqrt{2}/4  & 1/2        & \sqrt{2}/4 & 0           \\
  0           & \sqrt{2}/4 & 1/2        & \sqrt{2}/4  \\
  -\sqrt{2}/4 & 0          & \sqrt{2}/4 & 1/2
\end{array} \right)
\]

\begin{figure}[htp]\label{fig-octagonal}
\begin{center}
\includegraphics{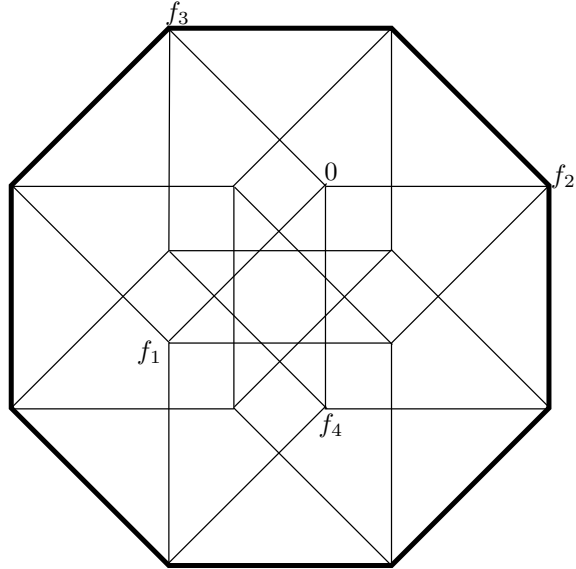}
\caption{The acceptance domain of the octagonal tiling, seen on $F$,
with the $1$-singular points.}
\end{center}
\end{figure}

Then the octagonal tiling is the tiling obtained from these data, with
canonical acceptance domain (that is, the projection of the $4$-dimensional
cube).
Figure~\ref{fig-octagonal} represents the acceptance domain
$K \subseteq F$ of this tiling, with the $1$-singular points.

We note $f_1, \ldots, f_4$ the projections of the $\{e_i\}$ on $F$.
The four singular hyperplanes are $H_i = \vect (f_i)$.
From the form of the matrix $\mat_{\mathrm{can}}{\pi_F}$, we see
that the $H_i$ have non-trivial stabilizers:
\[
\begin{array}{c}
\Stab_\Gamma (H_1) = \langle f_1, f_2 - f_4 \rangle \\
\Stab_\Gamma (H_2) = \langle f_2, f_1 + f_3 \rangle \\
\Stab_\Gamma (H_3) = \langle f_3, f_2 + f_4 \rangle \\
\Stab_\Gamma (H_4) = \langle f_4, f_1 - f_3 \rangle
\end{array}
\]
So we have for all $i \in \{1, \ldots, 4 \}$, $\alpha_i = 2$, so that
$\alpha = 2$, and the complexity of the octagonal tiling grows like $n^2$.

\subsection{Cubic Billiards}
The complexity of cubic billiards was already thoroughly studied by
mathematicians as well as computer scientists.
We show that our result gives the (well known) order of growth for the
complexity of a word obtained by a cubic billiard, even though we have no
control on the constants involved.

The assumption that $\Gamma := \pi_F (\Z^N)$ is dense in $F$, and $F$ is of
dimension $N-1$ implies that $\Gamma$ is of rank $N$.
The singular hyperplanes are $H_i := \vect\big( \{\pi_F(e_k)\}_{k \neq i}\big)$.
The rank of $\Stab_\Gamma (H_i)$ is at least $N-1$ (by definition), and
cannot be $N$, because it would imply $\Gamma \subseteq H_i$, and $\Gamma$
would not be dense in $F$.

So for all $i$, we have $\alpha_i = N-1$, and the complexity of a billiard
sequence in an $N$-dimensional hypercube in a suitable direction (which
ensures that $\Gamma$ is dense in $F$), grows like $n^{N-1}$.

This result was already known, with greater precision. See~\cite{Bar}.

\subsection{Discrete Hyperplane}
In the case $d = N-1$ (the Sturmian sequence being a particular case),
the space $F$ is one-dimensional, and there exists only one singular
hyperplane (the point $\{ 0 \} =: H_0$).

Then $\Stab_\Gamma (H_i) = \{ 0 \}$, and the complexity only depends of the
rank of $\Gamma$.
Let us assume $E \cap \Z^N$ is of rank $k$, that is the group of periods of
the resulting Delone set is of rank $k$.
Then $\Gamma$ is of rank $N-k$, and $\alpha_i = \rk \Gamma - 1 = N-k-1 = d-k$.

We deduce that the complexity of a Delone set of dimension $d$, periodic in
$k$ independent directions, obtained by cut and projection in a total space
of dimension $N = d+1$, grows like $n^{d-k}$.

In the case $k=0$ (no period), the complexity grows like $n^d$. This result is
in line with the results of rectangular complexity for discrete planes
of~\cite{BV}.

\section{Links Between Complexity and Cohomology}\label{sec:links}

\subsection{Cohomology for Tiling Spaces}

Cohomology first originates from algebraic topology. The idea of cohomology is
to associate Abelian groups (the \emph{cohomology groups})
to a topological space $X$. These groups are noted $(H^i(X))_{i \geq 0}$.
The construction is functorial (and contravariant) in the sense that given
a continuous morphism $f:X \rightarrow Y$, then there is an induced map
$f^* : H^i (Y) \rightarrow H^i (X)$.
Cohomology groups can be defined by various means on topological spaces
such as manifolds. We can mention for example the singular cohomology, the
De~Rham cohomology, of the \v{C}ech cohomology. For manifolds, these definitions
agree.
However, for spaces with less regularity, some definitions are more suited than
others.

The strength of cohomology is that it is a topological invariant. Homeomorphic
(or even homotopic) spaces have the same cohomology groups.
Furthermore, an exact sequence machinery exists which allows, to some extend,
to compute these groups actually.

Tiling spaces are quite complicated topological spaces: they are everywhere
locally the product of a Cantor set by an Euclidean ball.
It appeared that the good notion of cohomology for these spaces was \v{C}ech
cohomology. A remarkable feature of \v{C}ech cohomology is that it is well
behaved towards inverse limits: the $i$-th \v{C}ech cohomology group of an
inverse limit of spaces is the direct limit of the $i$-th cohomology groups of
these spaces.
This fact led to an actual computation of cohomology groups for a subclass of
tilings including some important examples. Anderson and Putnam~(\cite{AP})
showed that substitution tiling spaces were inverse limits of simple
complexes, which allowed to compute the associated invariants.

In the framework of cut and projection tiling spaces, another method exists and
was carried out by Forrest, Hunton and Kellendonk~(\cite{FHK}). This method uses
the fact that a group acts on the tiling space.
The cohomology groups actually computed with this method are the cohomology
groups of the group $\Z^N$ with value in a certain module of functions, and
this other notion of cohomology agrees with the \v{C}ech cohomology of the
topological space.
Their computations provided a criterion for infinitely generated cohomology,
which we will give below.

\subsection{Preliminaries and Statement of the Theorem}

We may wonder if complexity has a topological interpretation. Indeed,
intuitively, the more a tiling is complex, the less it should be repetitive.
As the usual distance which defines the topology of tiling spaces (see for
example~\cite{AP}) is linked to the notion of recurrence of patches, complexity
should have a dynamical or topological interpretation in terms of the associated
tiling space.
We will show that, for cut and projection tilings with almost-canonical
acceptance domain, one has such a link: the cohomology of a cut and projection
tiling space is finitely generated (over
$\Q$) if and only if its complexity grows as slow as possible.
More precisely, the statement is the following:
\begin{theo}\label{theorem-cohom}
Consider a cut and projection method $(\R^N = E \oplus F, \Z^N, \pi, K)$
which satisfies our assumptions.
Call $\Gamma = \pi_F (\Z^N)$.
Let $p$ be the associated complexity function.
Then there exists a constant $C_1 > 0$ such that $p(n) \geq C_1 n^\alpha$ with
$\alpha = d - (N-\rk\Gamma)$.
The cohomology groups over $\Q$ of the associated tiling space are finitely
generated if and only if the complexity is minimal, that is there exists $C_2$
such that $p(n) \leq C_2 n^\alpha$.
\end{theo}

In the following sections, we will see that this theorem is specific to cut
and projection Delone sets: such an equivalence is actually false in greater
generality, and we will give a counter example of one of these implications in
dimension one.

\begin{definition}
Given a $(N,d)$-cut and projection method, we call $\sg$ the set of singular
points, and $\sg_n$ the set of $n$-singular points.
These have a structure of cell complex, and we note $L^{(k)}$ the set of
$k$-cells of $\sg$ (for $0 \leq k \leq N-d-1$), and $L^{(k)}_n$ the set of
$k$-cells of $\sg_n$.
Remark that $L^{(k)}$ carries naturally a $\Gamma$-action, induced by the
translation.
\end{definition}

\begin{definition}\label{sing-k-space}
We call \emph{singular $k$-space} an affine space of dimension $k$, which is
singular, and is parallel to a space of the form $A = \bigcap_{j \in J}{H_j}$.
\end{definition}

\begin{theo}[\cite{FHK}]
The cohomology groups associated to a cut and projection method are finitely
generated if and only if $L^{(0)} / \Gamma$ is finite.
\end{theo}

The theorem~\ref{theorem-cohom} falls actually in two parts: first, it gives a
lower bound for the exponent of the complexity.
Then, it makes a link with cohomology. The cohomology groups are finitely
generated if and only if the exponent of the complexity is minimal.
Let us make a connection with our result on the minimal complexity exponent and
a conjecture by Lagarias and Pleasant.
\proc{Conjecture (\cite{LP}).}
Let $d$ the dimension be given. For each $1 \leq j \leq d$, there is a positive
constant $c_j = c_j(d,r,R)$ such that any Delone set $\Dron$ in $\R^d$ with
constants $(r,R)$ that satisfies
\[
 p(n) / n^{d+1-j} < c_j
\]
for all $n$ big enough has at least $j$ linearly independent periods.
\medbreak

What theorem~\ref{theorem-cohom} proves is that \emph{for a cut and projection
Delone set}, if $p(n) / n^{d+1-j}$ is bounded, then $d+1-j \geq d+\rk\Gamma-N$
which means $\rk\Gamma \leq N - j + 1$, which means that the Delone set has at
least $j-1$ independent periods.
The Delone set has at least $j$ independent periods if and only if the quantity
$p(n) / n^{d+1-j}$ tends to zero.
Our method gives a rather asymptotic result, and we do not have control on
the constants involved.
Therefore, even if we do not fully prove the conjecture in the case of cut and
projection Delone sets, the theorems in this article are seem to confirm the
result.
We do not know if our methods can be adapted to prove the conjecture.

\subsection{Proof of the Theorem}

Let us now prove theorem~\ref{theorem-cohom}. First, remark the following.
\begin{rem}
An element $x$ is in $L^{(0)}$ if and only if there exist $N-d$
singular hyperplanes $(H_j)_{j \in J}$ the intersection of which is $\{0\}$,
elements $(\gamma_j)_{j \in J}$ of $\Gamma$, and elements $v_j \in V_j$, such
that:
\[
 x = \bigcap_{j \in J}{\big( H_j + \gamma_j + v_j \big)}
\]
\end{rem}

Consider a $(N,d)$-cut and projection method, with complexity function $p(n)$.
We proved (theorem~\ref{theorem}) that $p(n)$ grows like $n^\alpha$.
We want to prove that $\alpha \geq d + \rk \Gamma - N$, with equality if and
only if $L^{(0)} / \Gamma$ is finite.

\proc{Notation.}
Recall that $\Gamma^j$ is the stabilizer in $\Gamma$ of the singular
hyperplane $H_j$.
For $J \subset \{1 \ldots m\}$, we will note $\Gamma^J$ the stabilizer of the
intersection $\bigcap_{j \in J}{H_j}$, and $a_J = \rk{\Gamma^J}$.
We use the convention $\Gamma^\emptyset = \Gamma$.
\medbreak

\proc{Demonstration of Theorem~\ref{theorem-cohom}.}
We use the notations of the theorem~\ref{theorem}:
$\alpha = \sum_{i \in I}{\alpha_i}$, where
$\alpha_i = \rk{(\Gamma / \Gamma^i)} - 1 = \rk{\Gamma} - a_i - 1$, and the
$\alpha_i$, $i \in I$ are chosen to be the $N-d$ greatest amongst the
$\alpha_j$, $1 \leq j \leq m$, and such that $\bigcap_{i \in I}{H_i} = \{0\}$.
Therefore, we have:
\begin{equation}\label{eqn-alpha-min}
 \begin{array}{lcl}
  \alpha &=& (N-d) \rk{\Gamma} - (N-d) - \displaystyle\sum_{i \in I}{a_i} \\
         &=& \big( d + \rk{\Gamma} - N \big) + \left[ (N-d-1) \rk{\Gamma} - 
                   \displaystyle\sum_{i \in I}{a_i} \right]
 \end{array}
\end{equation}
and the $a_i$, $i \in I$, are chosen to be the $N-d$ \emph{smallest} amongst
the $a_j$, $1 \leq j \leq m$ (with the transversality condition).

Given two sets of indices $J$, and $J'$ in $\{0, \ldots, m \}$, consider the
following application:
\[
 \begin{array}{cccl}\displaystyle
  \phi_{J,J'} : & \Gamma / \Gamma^{J \cup J'} & \longrightarrow &
     \Gamma / \Gamma^J \times \Gamma / \Gamma^{J'} \\
          & [x]_{J \cup J'} & \longmapsto & ([x]_J , [x]_{J'})
 \end{array}
\]
where $[x]_J$ stands for the class of $x$ in $\Gamma$ modulo $\Gamma^J$.
One checks easily that this application is well defined and one-to-one.
So the rank of the target group is equal or greater than the rank of the
source group (remark that both of these groups are free Abelian).
Therefore, the following condition on the ranks holds:
\[
 \rk{\Gamma} - a_{J \cup J'} \leq \rk{\Gamma} - a_J + \rk{\Gamma} - a_{J'}
\]
which gives:
\begin{equation}\label{ineq-ranks}
a_J + a_{J'} - a_{J \cup J'} \leq \rk{\Gamma}
\end{equation}

Let us iterate this inequality on equation~(\ref{eqn-alpha-min}) in order to
deduce the first part of the theorem.
We label $i(1), \ldots, i(N-d)$ the elements of $I$.
Then we write:
\begin{equation}\label{iteration}
 \begin{array}{rcl}
  a_{i(1)} + a_{i(2)} & \leq & a_{\{i(1),i(2) \}} + \rk{\Gamma} ;\\
  a_{i(1)} + a_{i(2)} + a_{i(3)} & \leq & a_{\{i(1),i(2) \}} + a_{i(3)}
    + \rk{\Gamma} \\
      & \leq & a_{\{i(1), i(2), i(3)\}} + 2 \rk{\Gamma} ;\\
      & \vdots  & \\
  \displaystyle\sum_{i \in I}{a_i} & \leq & (N-d-1) \rk{\Gamma} + a_I \\
      & = & (N-d-1) \rk{\Gamma}.
 \end{array}
\end{equation}
We deduce that $\left[ (N-d-1) \rk{\Gamma} - \sum_{i \in I}{a_i} \right]$ is
positive, which yields:
\[
 \alpha \geq d + \rk{\Gamma} - N
\]
This proves the first part of theorem~\ref{theorem-cohom}.

Let us now prove that $\alpha = d-N+\rk\Gamma$ if and only if $L^{(0)} / \Gamma$
is finite.
First, we prove that the equality implies finiteness of $L^{(0)}/\Gamma$, by
a decreasing induction.
Assume $\alpha = d + \rk{\Gamma} - N$, and consider the following recurrence
hypothesis:
for a given $k < N-d$ and every space of the form $A = \bigcap_{j \in J}{H_j}$
of dimension $k$ ($\Card{J} = N-d-k$), then the set of singular $k$-spaces
(recall definition~\ref{sing-k-space}) which are directed by $A$ is the union of
only finitely many $\Gamma$-orbits of translates of $A$.

The assertion above is obviously true for $k = N-d-1$, since the set of singular
$(N-d-1)$-spaces directed for by $H_i$ is included in $H_i + V + \Gamma$ for
all $i$, and $V$ is finite.

Let us assume that the recurrence hypothesis is true for a given $k+1$, and let
us prove that the set of singular $k$-spaces directed by a given space of
dimension $k$ is the union of finitely many $\Gamma$-orbits.
Let $A = \bigcap_{i \in J}{H_i}$ be a vector subspace of dimension $k$, with
$\Card{J} = N-d-k$.
We choose $j \in J$, and define $J' := J \setminus j$, and
$A' = \bigcap_{j \in J'}{H_j}$.
Then consider the map $\phi_{J', j} := \phi_{J',\{j\}}$ as defined above.
% We assumed that $\alpha = d + \rk{\Gamma} - N$. Therefore, we have:

We deduce from the iterative reasoning of equation~(\ref{iteration}) that
for all set $I'$ of $N-d$ elements such that $\cap_{i \in I'}{H_i} = \{0\}$,
the following holds:
\[
 (N-d-1) \rk{\Gamma} - \sum_{i \in I'}{a_i} \geq 0
\]
(equation~(\ref{iteration}) is not specific to the choice of $I$).
Furthermore, we assumed that $\alpha = d + \rk\Gamma - N$. Therefore, for the
specific choice $I' = I$, the inequality above is an equality.
But the choice of $I$ is the choice for which the $a_i$ are \emph{minimal}.
Therefore, the above expression is positive, and its maximum is zero.
So for all choice of $I'$ such that $\cap_{i \in I'}{H_i} = \{0\}$:
\begin{equation}
\label{eqn-zero}
 (N-d-1) \rk{\Gamma} - \sum_{i \in I'}{a_i} = 0
\end{equation}

We claim that the following equality holds, which means
inequality~(\ref{ineq-ranks}) is in fact an equality:
\begin{equation}\label{eq-ranks}
 a_{J'} + a_j = a_J + \rk{\Gamma}
\end{equation}
for if we had $a_{J'} + a_j < a_J + \rk{\Gamma}$, then we could choose $I'$ to
be a set of $N-d$ elements of $\{1, \ldots, m\}$, containing $J$.
Then, applying the inequalities of~(\ref{iteration}) to $I'$, we would obtain a
strict inequality at some step, and so we would obtain:
\[
 (N-d-1) \rk{\Gamma} - \sum_{i \in I'}{a_i} > 0
\]
in contradiction equation~(\ref{eqn-zero}) above.

Consider the map $\phi_{J',j}: \Gamma / \Gamma^J \rightarrow (\Gamma /
\Gamma^{J'}) \times (\Gamma / \Gamma^j)$: it is one-to-one.
Furthermore, the source group has rank $\rk{\Gamma} - a_J$ and the target group
has rank $2 \rk{\Gamma} - a_{J'} - a_j$.
Those have the same rank, by equation~(\ref{eq-ranks}).
It doesn't necessarily mean that $\phi_{J',j}$ is onto, but its image is of
finite index, so its cokernel
$G := (\Gamma / \Gamma^{J'}) \times (\Gamma / \Gamma^j)$ is finite.

We consider the set of $k$-dimensional affine spaces of the following form:
\[
 (A' + \gamma_1) \cap (H_j + \gamma_2)
\]
with $\gamma_1, \gamma_2 \in \Gamma$.
Due to stabilization of the spaces involved, it is enough to consider
$\gamma_1 \in \Gamma^{J'}$ and $\gamma_2 \in \Gamma^j$.
Then the set of translates of singular $k$-spaces of the form above is in
correspondence with $(\Gamma/\Gamma^{J'}) \times (\Gamma/\Gamma^j)$.
Quotiening by the action of $\Gamma$ (or by the image of $\phi_{J',j}$, which
amounts to the same), we have that the set of $\Gamma$-orbits of $k$-dimensional
spaces of the form above is finite.

In general, the singular $k$-spaces parallel to $A$ are of the form:
\[
 (A' + \gamma_1 + v_1) \cap (H_j + \gamma_2 + v_2)
\]
where $A'$ is chosen amongst the (finite) set of $(k+1)$-dimensional vector
spaces containing $A$ and of the form $\bigcap_i{H_i}$, $j$ is chosen
amongst the hyperplanes such that $H_j \cap A' = A$,
$\gamma_1 \in \Gamma /\Gamma^{J'}$, $\gamma_2 \in \Gamma / \Gamma^j$, $v_2$ is
chosen in $V$ (the set of vertices of the acceptance domain), which is finite.
Plus, by recurrence hypothesis, there are only finitely many $\Gamma$-orbits
of translates of singular $k+1$-spaces directed by $A'$, so $v_1$ can be chosen
in a finite set.
Thus, as all choices we made were finite, the number of $\Gamma$-orbits of
singular $k$-spaces directed by $A$ is finite.
We can say this for every singular $k$-space $A$, so the recurrence hypothesis
holds for $k$.

We conclude, using induction down to $k=0$, that $L^{(0)}$ is an union of
finitely many $\Gamma$-orbits.

Conversely, assume that $\alpha > d - N + \rk\Gamma$.
It means that there is at least one inequality in equations~(\ref{iteration})
which is not an equality.
It means that there is $J' \subseteq \{1, \ldots, m\}$ and
$j \in \{1, \ldots, m\}$ such that the $(H_i)_{i \in J\cup\{j\}}$ are
transverse. We note, as above, $J = J' \cup \{j\}$. Then,
\[
 a_{J'} + a_j  < a_J + \rk\Gamma
\]
As we said above, it proves that the set of singular $k$-spaces of the form:
\[
 \left( \bigcap_{i \in J'}{H_i} + \gamma_1 \right) \cap (H_j + \gamma_2)
\]
is the union of infinitely many $\Gamma$-orbits of translates of
$\bigcap_{j \in J}{H_j}$.
Therefore, taking the intersection of these spaces with any complementary
singular subspace, we deduce that $L^{(0)}$ is made of infinitely many
$\Gamma$-orbits.
\ep\medbreak

\section{The case of Dimension One}

In dimension one, for general tilings (not necessarily obtained by cut and
projection), the links between cohomology and complexity are quite different.

To a tiling of $\R$, one can associate a bi-infinite word $w$, which
corresponds to the sequence of its tiles.
In this section $p(n)$ stands for the number of subwords of length $n$ of
a given word $w$.
Beware that before this section, $p(n)$ in the one dimensional case was defined
as the number of words of length $2n$.
These two definitions give of course the same order of growth, only the
constants are changed.

\subsection{Definitions and construction}
Let us fix the notations we will use in the following sections.
We will consider a finite alphabet $\Aron$.
Given $w \in \Aron^\Z$ a bi-infinite word over $\Aron$, we write:
\[
 w = \ldots w_{-2} w_{-1} . w_0 w_1 w_2 \ldots
\]
The space $\Aron^\Z$ is a topological space in a natural way (it is a Cantor
set). The group $\Z$ acts naturally on $\Aron^\Z$ by the \emph{shift}
homeomorphism. This homeomorphism, $\tau$, is defined as follows:
\[
\tau (\ldots w_{-2} w_{-1} . w_0 w_1 \ldots) =
   (\ldots w_{-1}  w_0 . w_1 w_2 \ldots)
\]

\begin{definition}
 We call \emph{subshift} of $\Aron^\Z$ generated by $w$, the closure in
$\Aron^\Z$ of the set:
\[
 \Theta_w := \{ \tau^n (w) \tq n \in \Z \}
\]
together with the shift application $\tau$.
\end{definition}

Such a subshift is totally disconnected, compact, and $\tau$-stable.
We now define the suspension of a subshift $\Theta_w$.
\begin{definition}
If $(\Theta_w, \tau)$ is the subshift of $\Aron^\Z$ associated to $w$, we
define its \emph{suspension} $(\Omega_w, \phi)$ by:
\[
\Omega_w = \Theta_w \times \R / \sim
\]
where $\sim$ is the equivalence relation generated by the relations
$(x,t) \sim (\tau(x),t-1)$ for all $x \in \Theta_w$ and all $t \in \R$.
The map $\phi$ is an action of $\R$ defined by:
\[
\phi_s ( [(x,t)] ) := [(x,t+s)]
\]
With these definitions, $(\Omega_w, \phi)$ is a $\R$-dynamical system.
\end{definition}

\begin{rem}
The construction of the suspension allows us to make the link between words
and tilings of $\R$. An element $[(x,t)]$ of $\Omega_w$ can be seen as a tiling
of $\R$ by collared tiles, which would all be of length one, the set of colors
being $\Aron$. The tiling follows then the pattern of the word $x$ and the
parameter $t$ indicates the position of the tiling.
There is a natural action of $\R$ by translation on tilings: this action
corresponds to $\phi$ on the suspension space.
\end{rem}

Again, $\Omega_w$ has a natural topology, which we define now.
\begin{definition}
The following subsets of $\Omega_w$ define a basis of open sets for the
topology:
\[
 U(V,I) = \{ [(x,t)] \tq x \in V, \  t \in I \}
\]
with $V$ an open subset of $\Theta_w$ and $I$ an open interval of $\R$.
\end{definition}

The following proposition is classical, and we state it without
demonstration.
\begin{prop}
With this topology, $\Omega_w$ is compact and connected.
\end{prop}

Now, let us define the Rauzy graphs associated to $w$.
\begin{definition}
For all $n \in \N$, define $F_n(w)$ the set of factors of $w$, of size $n$,
that is the set of all finite subwords of $w$ of length $n$.
The $n$-th Rauzy graph $R_n$ associated to $w$ is an oriented graph defined by:
\begin{itemize}
\item a set of vertices $V_n = F_n (w)$;
\item a set of oriented edges $E_n$: there is an edge from
$a_1 \ldots a_n \in V_n$ to $b_1 \ldots b_n \in V_n$ if $a_i = b_{i-1}$ for
$2 \leq i \leq n$, and if the word $a_1 \ldots a_n b_n$ is a factor of $w$.
\end{itemize}
We note $_a f_b$ the oriented edge from $af$ to $fb$ where $a,b \in \Aron$
and $f \in F_{n-1}$, and such that $af, fb \in V_n$. Such an edge exists if and
only if $afb \in F_{n+1}$.
\end{definition}

Rauzy graphs are very well suited to study combinatorial properties on
words. Namely, the vertices with more than one outgoing edge correspond exactly
to the factors which can be extended to the right in more than one way (such
factors are sometimes called ``special factors'').
As all factors of $w$ can be extended to the right or to the left, the Rauzy
graphs $R_n$ are connected by oriented paths (or strongly connected).

Remark that these graphs are combinatorial objects. However, we can (and will)
regard them as metric spaces as well: a graph is the disjoint union indexed
over the edges of copies of the interval $[0,1]$, quotiented by the adjacency
condition (gluing the vertices together).

\subsection{Inverse Limits of Rauzy Graphs}

We will first prove that one-dimensional tiling spaces are inverse limits of
Rauzy graphs.
It is already known that they are inverse limits of graphs, namely the G\"ahler
complexes.
Even though Rauzy graphs and G\"ahler complexes are constructed in a very
similar way, they can be different.
Rauzy graphs are best suited to capture combinatorial information on the word,
such as complexity or the number of ``bi-extendable'' words (or special
factors).
Therefore, they are a convenient technical tool if we want to make a link
between complexity and cohomology.

In order to build inverse limits, we need projection maps
$R_{n+1} \rightarrow R_n$.
The simple fact that a word of length $n+1$ can be shortened to a word of
length $n$ will allow us to define maps $\gamma_n : R_{n+1} \rightarrow R_n$.
In this construction, one could shorten the word from the left or from the
right. We make the following choice:
\begin{definition}
 Define $\gamma_n : R_{n+1} \rightarrow R_n$ in the following way.
If $n$ is even, $\gamma_n$ is defined by:
\begin{itemize}
\item $\gamma_n \big( a := a_1 \ldots a_{n+1} \big) = a_1 \ldots a_n$ if
$a \in S_{n+1}$ (remove a letter on the right).
\item If an edge $e \in E_n$ goes from $a$ to $b$, then there exists in $R_n$
an edge from $\gamma_n ({a})$ to $\gamma_n ({b})$.
The application $\gamma_n$ then maps $e$ on that edge, and its restriction on
$e$ is the identity map on the segment $[0,1]$.
\end{itemize}
If $n$ is odd, $\gamma_n$ is defined in a similar way by:
\begin{itemize}
\item $\gamma_n \big( a := (a_1 \ldots a_{n+1}) \big) = a_2 \ldots a_{n+1}$, if
${a} \in S_{n+1}$ (remove a letter on the left).
\item The definition of $\gamma_n$ on the edges is exactly the same.
\end{itemize}
\end{definition}

\begin{prop}
 The maps $\gamma_n$ are continuous.
\end{prop}

\proc{Demonstration.}
Straightforward.
\ep\medbreak

We can now consider inverse limits of the Rauzy graphs, with respect to the
maps $\gamma_n$.
\begin{definition}
 The inverse (or projective) limit of the complex $(R_n, \gamma_n)_{n \in \N}$
is the following subset of $\prod_{n \in \N}{R_n}$:
\[
 \varprojlim{\big( R_n, \gamma_n \big)_{n \in \N}} :=
\left\{ (x_1, x_2, \ldots) \in \prod_{n \in \N}{R_n} \tq \forall n \in \N, \ 
\gamma_n (x_{n+1}) = x_n \right\}
\]
This inverse limit is a topological space with the relative topology induced
by the product topology of $\prod_{n \in \N}{R_n}$.
It is a compact set, being a closed subset of a compact set.
\end{definition}

\begin{theo}\label{inv-lim}
Let $w$ be a bi-infinite word over the finite alphabet $\Aron$. Let $\Theta_w$
be its associated subshift of $\Aron^\Z$.
Note $\Omega_w$ the suspension of $(\Theta_w, \tau)$. We then have an
homeomorphism:
\[
 \Omega_w \simeq \varprojlim{(R_n, \gamma_n)_{n \in \N}}
\]
\end{theo}

For the proof, we use the now classical methods developed in~\cite{AP}.

\proc{Demonstration.}
Let us build an homeomorphism between $X := \varprojlim_n{R_n}$ and
$\Omega := \Omega_w$.
For all $n \in \N$, we will define an application $\psi_n : \Omega \rightarrow
R_n$. It will allow us to define an application $\psi : \Omega \rightarrow X$,
which we will prove, is a homeomorphism.
Let $P \in \Omega$. Then $P = (m,t)$ with $m \in \Theta_w$, and $t$ can be
chosen in $\left[0,1\right[$.
Note $m = (m_i)_{i \in \Z}$.
Define $\psi_n (P)$ as follows:
\begin{itemize}
 \item if $n$ is even ($n = 2p$), consider the vertices $s_1$ and $s_2$ of
$R_n$ defined by the following words:
\[
 \begin{array}{lcl}
  s_1 & := & m_{-p} \ldots m_{p-1} \\
  s_2 & := & m_{-p+1} \ldots m_p
 \end{array}
\]
There exists an oriented edge isometric to $[0,1]$ from $s_1$ to $s_2$.
Then $\psi_n (P)$ is defined to be the point in position $t \in [0,1]$ in this
edge.
 \item if $n$ is odd, ($n = 2p+1$), consider the vertices $s_1$ and $s_2$ of
$R_n$ defined by the following words:
\[
 \begin{array}{lcl}
  s_1 & := & m_{-p-1} \ldots m_{p-1} \\
  s_2 & := & m_{-p} \ldots m_{p}
 \end{array}
\]
$\psi_n (P)$ is defined to be the point in position $t \in [0,1]$ in the edge
from $s_1$ to $s_2$.
\end{itemize}

A direct computation shows that the $\psi_n$ are compatible with the
maps $\gamma_n$, that is for all $n \in \N$:
\[
 \psi_n = \gamma_{n+1} \circ \psi_{n+1}
\]

Then, due to this compatibility relations, the application $\psi_n$ define a
map $\psi : \Omega \rightarrow X$.
Let us check that $\psi$ is one-to-one, onto, and continuous. By compactness of
$\Omega$, it will be a homeomorphism.

Let $P, P' \in \Omega$. We write $P = (m,t)$ and $P' = (m',t')$, with $t, t'
\in \left[0,1\right[$.
We assume that $P \neq P'$.
If $t \neq t'$, it is clear that $\psi_n (P) \neq \psi_n (P')$ for all $n$, and
hence $\psi (P) \neq \psi (P')$.
If $t = t'$, then $m \neq m'$, and hence there is some $p$ such that
\[
 m_{-p-1} \ldots m_{p-1} \neq m'_{-p-1} \ldots m'_{p-1}
\]
Hence, $\psi_{2p+1}(P)$ and $\psi_{2p+1}(P')$ do not belong to the same edge of
$R_{2p+1}$, and so they cannot be equal.
In all cases, if $P \neq P'$, then $\psi(P) \neq \psi (P')$, and $\psi$ is
one-to-one.

Let us prove that $\psi$ is onto. Given $(x_n)_{n \in \N}$ an element of
$X$, let us construct a preimage in $\Omega$.
By construction, there exists at most one edge from one given vertex of a Rauzy
graph to another.
Therefore, the element $x_n \in R_n$ is entirely determined by two vertices
$s_n$ (starting point), $s'_n$ (ending point), and an element
$t_n \in \left[0,1\right[$ (position on the edge).
It is obvious that $t_n =: t$ does not depend on $n$ (because of the
compatibility conditions).
From the sequence of the $s_n$, we will construct a bi-infinite word:
by compatibility condition, each $s_n$ extends $s_{n-1}$.
We define $m$ to be the bi-infinite word obtained by increasing union of these
finite words. The factor $s_n = s_n^{(1)} \ldots s_n^{(n)}$ will be found in the
following position (say, for $n$ even):
\[
 m = \ldots s_n^{(1)} \ldots s_n^{(n/2)} . s_n^{(n/2 + 1)} \ldots s_n^{(n)}
\ldots
\]
where the dot lies just before the zero-th letter of $m$.
Then we check that $(m,t) \in \Omega$ is a preimage for $(x_n)_{n \in \N}$.

Finally, let us prove that $\psi$ is continuous.
Let $P \in \Omega$, and $V$ be an open neighborhood around $x := \psi(P)$.
Then $V$ is of the form:
\[
 V = V_1 \times \ldots \times V_n \times \prod_{i>n}{R_i}
\]
where the $V_i$ are neighborhoods of $x_i := \psi_i (P)$ in $R_i$.
Let us prove that $\psi^{-1}(V)$ contains an open set containing $P$.
There are two cases to consider:
\begin{itemize}
 \item If $P = [(m,t)]$ with $t \in \left]0 , 1 \right[$, then we can restrict
all the $V_i$ such that for all $i \leq n$, $V_i$ is an open neighborhood of
$x_i$ which is entirely included in an edge of $R_n$.
Furthermore, if we note, as above, $x_i = (s_i, s'_i, t)$, the union
of $s_n$ and $s'_n$ define a finite word of length $n+1$ centered in zero,
which, in turn, defines a clopen set of $\Theta_w$ (that is a set which is both
open and closed). We call $W$ this clopen set.
Finally, as the $V_i$ are all included on edges isometric to $[0,1]$, each of
these neighborhoods define an interval $I_i$ of $\left]0 , 1 \right[$.
We call $I$ their (finite) intersection, which is an open interval containing
$t$.
Then $U(W,I)$ is an open set of $\Omega$ included in $\psi^{-1}(V)$, and
containing $P$.
 \item If $P = [(m,0)]$, the construction is about the same. We construct a
finite subword of $m$, which defines a neighborhood $W$ of $\Theta_w$, and such
that $\psi \big( [(W,0)] \big) \subseteq V$.
Then, we remark that the $V_i$'s define (up to isometry), a certain number of
neighborhoods of $1$ in $[0,1]$ which correspond to the incoming edges, and a
certain number of neighborhoods of $0$ in $[0,1]$ which correspond to the
outgoing edges.
We call $I_k$ all the neighborhoods of $1$ (when $1 \leq i \leq n$), and $J_k$
the neighborhoods of $0$.
Now, define $I$ the intersection of all the $I_k$, and $J$ the intersection of
all the $J_k$.
Finally, call $I' = (J-1) \cup I$, which is an open neighborhood of $0$ in $\R$.
Then $U(W,I')$ is an open set of $\Omega$ included in $\psi^{-1}(V)$, and
containing $P$.
\end{itemize}
We proved that $\psi$ is continuous. Hence, by compacity of $\Omega_w$,
it is an homeomorphism, and the theorem is proved.
\ep\medbreak

\subsection{Preliminaries}\label{preliminaries}

We proved (theorem~\ref{inv-lim}) that tilings spaces can be seen as inverse
limits of some nice topological spaces, namely the Rauzy graphs, under the maps
$(\gamma_n)_{n \in \N}$.

The number of vertices of the $n$-th Rauzy graph $R_n$ equals the complexity
$p(n)$, and its number of edges is related to the function $p(n+1)-p(n)$.
This makes it suited to study complexity.
In comparison, the number of edges of $G_n$, the $n$-th G\"ahler complex equals
$p(n)$, but we have no real control on the number of vertices: this makes it
less
suited for our purpose.

Let us make actually the link between complexity and the combinatoric properties
of the Rauzy graphs associated to a word.
As in the sections above, $w$ is a given bi-infinite word over a finite
alphabet $\Aron$.
Remember that $p(n)$ counts the number of different factors of $n$ letters
which appear in $w$.

\begin{definition}
We call $s(n) = p(n+1) - p(n)$.
\end{definition}

Remark that a factor of $w$ of length $n+1$ can be obtained from a factor of
length $n$ by extending it on the right by a suitable letter. Any factor of
$F_{n+1}$ can be obtained in this way.
This simple fact gives the following proposition.
\begin{prop}
\[
 s(n) = \sum_{w \in F_n}{[\textrm{number of ways to extend } w \textrm{ on the
right } - 1]}
\]
\end{prop}
So the function $s$ counts the number of ``exceptional'' ways to extend words on
a given side, and in a way that the extension is still a subword of $w$.

\proc{remark}
Over a two-letter alphabet, $s(n)$ simply counts the number of ``special
factors'' of length $n$ of $w$: all subword of $w$ of a given length can be
extended in at least one way. The special factors are those which can be
extended in more than one way.
\medbreak

Let us now introduce a notation which will make easier the formulation of the
notion of branching point.
\proc{Notation.}
Let $G$ be an oriented graph.
 \begin{itemize}
  \item If $v$ is a vertex of $G$ , we note $d^+ (v)$ the
   number of outgoing edges starting from $v$.
   We define $d^- (v)$ similarly as the number of incoming edges.
  \item If $e$ is an edge of $G$, we call $e^+$ the vertex
   towards which $e$ is directed (so that $e$ is an incoming edge with respect
   to the vertex $e^+$).
 \end{itemize}
\medbreak

We now make explicit the link between $s(n)$ and the cohomology of the
$n$-th Rauzy graph $R_n$.
\begin{prop}\label{bound-cohom-b}
The rank of the cohomology $H^1(R_n;\Z)$ of the $n$-th Rauzy graph equals
$s(n)+1$.
\end{prop}

\proc{Demonstration.}
Let $n \in \N$. We call $V$ the set of vertices of $R_n$ and $E$ the set of its
edges.
By construction $\Card{V} = p(n)$.
Every factor of $w$ of length $n+1$ can be obtained by extending by one letter
to the right a word of length $n$. That is, by choosing a vertex $v \in V$ and
an outgoing edge from $v$.
Plus, every such choice of vertex and edge provides a different subword of $w$
(by construction of the Rauzy graphs).
Therefore, we can write:
\[
 p(n+1) = \sum_{v \in V}{d^+(v)} = \Card{E}
\]
Furthermore, it is classical for connected graphs that:
\[
 \rk{H^1(R_n;\Z)} = \Card{E} - \Card{V} + 1
\]
We deduce:
\[
  \rk{H^1 (R_n ; \Z)} = p(n+1) - p(n) + 1 = s(n) + 1
\]
\ep\medbreak

We now state a key lemma which links cohomology in the Rauzy graphs and in
their inverse limit.
\begin{lemma}\label{lemma-cohom-ccl}
We assume that there exists $M \in \N$ such that the rank of the cohomology
of infinitely many Rauzy graphs of the word $w$ is bounded by $M$.
Then the \v{C}ech cohomology over $\Q$ of the suspension space of $w$ is
generated by at most $M$ elements.
\end{lemma}

\proc{Demonstration.}
We know that the suspension space of $w$ is homeomorphic to the inverse limit of
the Rauzy graphs.
It implies that the \v{C}ech cohomology of the suspension space is the direct
limit of the cohomology of the graphs.
Let $x_1, \ldots, x_{M+1}$ be $M+1$ elements in the cohomology of the limit.
We call $\hat{x}_1, \ldots, \hat{x}_{M+1}$ lifts of the $x_i$, so that
$\hat{x}_i$ is an element of the cohomology of the complex $R_{\phi(i)}$.
Let $N$ be an upper bound for all the $\phi(i)$.
We can assume that $N$ is chosen in such a way that the cohomology of $R_N$
is generated by less than $M$ elements.
We call $y_1, \ldots, y_{M+1}$ the images of the $\hat{x}_i$ in the cohomology
of the graph $R_N$.
The rational cohomology of $R_N$ is a vector space of dimension at most $M$.
Therefore, there exists a linear relation between the $(y_i)_{i=1}^{M+1}$.
Taking the image of this linear relation in the rational cohomology of the
suspension space, we deduce that the $x_i$ do not form a free family.
\ep\medbreak

\subsection{Simple Words Have Finitely Generated Cohomology}

We will now prove that if a word is not too complex (that is if its
complexity function $p(n)$ is bounded by $Cn$ with $C$ a constant), then
the cohomology of its associated subshift is rationally finitely generated.

\begin{lemma}\label{lemma-bn}
Let $w$ be a word over a finite alphabet, and $p$ the corresponding complexity
function.
We assume that there exists a constant $C$ such that for all $n$,
\[
p(n) \leq C n
\]
Then the function $s$ does not tend to infinity.
\end{lemma}

Actually, a much stronger result was proved by Cassaigne: the function
$s$ is bounded (under some assumptions on the word, see~\cite{Cas} for this
result).
However, we just need this weaker result to prove
proposition~\ref{simpl-finite}.

\proc{Demonstration.}
It is enough to remark that $p(n+1) = p(n) + s(n)$, so that
$p(n) = \sum_{k=0}^{n-1}{s(k)}$.
Assume $s(k)$ is greater than $M+1$ for all $k$ big enough. Then $p(n)$ is
greater than $Mn$ for all $n$ big enough.
If $s(k)$ tends to infinity, $M$ can be chosen as big as one wants, which
implies that $p$ is not dominated by a linear function.
\ep\medbreak

With this, theorem~\ref{inv-lim}, and the results of
section~\ref{preliminaries}, we can prove the following proposition.

\begin{prop}\label{simpl-finite}
Let $w$ be a bi-infinite word over a finite alphabet.
Let $p$ be the complexity function associated to $w$. We assume that
$p(n) \leq Cn$ for all $n$ big enough.
Then the first rational cohomology group of the associated suspension space
is finitely generated.
\end{prop}

\proc{Demonstration.}
By theorem~\ref{inv-lim}, the suspension space of $w$ is the inverse limit of
its Rauzy graphs.
It implies that its first cohomology group is the direct limit of the
first cohomology groups of the complexes.
However, $s(n)$ doesn't tend to infinity, by lemma~\ref{lemma-bn}.
Therefore, by proposition~\ref{bound-cohom-b}, there exists an $M \in \N$
such that the first cohomology group of infinitely many Rauzy graphs
is generated by less than $M$ elements.
Finally, lemma~\ref{lemma-cohom-ccl} allows us to conclude that the first
cohomology group of the suspension space associated to $w$ is finitely
generated (over $\Q$).
\ep\medbreak

\subsection{Complex Tilings Can Have an Homologicaly Simple Tiling Space}

\proc{Notation.}
Let $f$ be a finite word over the finite alphabet $\Aron$.
We note $\abs{f}$ the length of $f$.
\medbreak

We now want to give an example of minimal, recurrent word, which has an
exponential complexity, and such that its associated tiling space has finitely
generated cohomology over $\Q$.
The idea is to find a word of high complexity, but such that infinitely many
of its Rauzy graphs are homotopic to a wedge of two circles, and
conclude thanks to lemma~\ref{lemma-cohom-ccl}.

This example was hinted to me by Julien Cassaigne, during a discussion we had
when he came to Lyon.

Let $\{a,b\}$ be an alphabet of two letters.
We choose a finite word in which all factors of length $2$ over $\{a,b\}$
appear, except $aa$.
We call $f_1$ this word.
We define $u_1 = a a b f_1 b a a$ and $v_1 = a b b f_1 b a a$.

Assuming we constructed an $f_k$, $v_k$ and $u_k$ of same %odd
length, let us construct $f_{k+1}$, $v_{k+1}$ and $u_{k+1}$ by induction.
We define a word $f_{k+1}$ over the alphabet $\{0,1\}$ in the following way:
$f_{k+1}$ is a finite word were all factors of length $\abs{u_k}$ without
$00$ over the alphabet $\{0,1\}$ appear.
Then, apply the following substitution on $f_{k+1}$:
\begin{align*}
 0 & \mapsto u_k \\
 1 & \mapsto v_k
\end{align*}
We still call $f_{k+1}$ the word we obtained after substitution. It is now a
word over $\{a,b\}$.
Then, we define:
\begin{align*}
u_{k+1} & = u_k u_k v_k f_k v_k u_k u_k \\
v_{k+1} & = u_k v_k v_k f_k v_k u_k u_k
\end{align*}
Then, $u_{k+1}$ and $v_{k+1}$ have the same length.

For all $k$, the word $u_k$ appears in $f_k$ (at various places), and so it
appears in the middle of $u_{k+1}$.
By making choices of the position of $u_k$ in $f_k$ for all $k$, we can form a
sequence of finite words over $\{a,b\}$, increasing on both ends.
Such a sequence is not unique, as $u_k$ appears many times in $f_k$, and so many
choices can be done.
We can take a limit to obtain a bi-infinite word over $\{a,b\}$.
We choose a word given by this method, and we call it $w$.

As an aside, we remark that one can prove that $w$ is uniformly repetitive, that
is, for all finite factors $f$ of $w$, there exists a $N \in \N$ such that an
occurrence of $f$ appears in every subword of length $N$ of $w$.

\begin{prop}
Let $p$ be the complexity function of $w$.
Then $p(n)$ is not dominated by any polynomial function in $n$.
\end{prop}

\proc{Demonstration.}
Let $p'$ be the complexity function of the subshift of $\{a,b\}$ without $aa$.
This subshifts contains $\{bab,bbb\}^\Z$, so we have $p'(3n) \geq 2^n$, which
means that $p'$ is exponential.

Let us now compare $p$ and $p'$.
By definition, $p(\abs{u_n}^2) \geq p'(\abs{u_n}) \geq 2^{\abs{u_n}/3}$.
Therefore, for some values of $N \in \N$, $p(N) \geq 2^{\sqrt{N}/3}$.
It proves that infinitely often, $p$ is greater than any polynomial function.
\ep\medbreak

Now, let us prove that the Rauzy graphs of $w$ are infinitely often
homotopic to a wedge of two circles.
\begin{lemma}
For infinitely many $N \in \N$, there is only one right special factor of length
$N$ in $w$ (\ie a factor which can be extended in more than one way on the
right).
\end{lemma}

\proc{Demonstration.}
Remark that $w$ has a very specific structure: one can uniquely rewrite $w$ as a
word over the alphabet $\{ u_1, v_1 \}$.
To do so, mark each occurrence of $aaaa$ or $baaab$, and cut it by adding a
dot:
$aaaa \rightarrow aa.aa$ and $baaab \rightarrow baa.ab$.
We then see $w$ as a concatenation of blocks, these blocks being precisely
$u_1$ and $v_1$.
By this method, we can rewrite $w$ as a word over the alphabet $\{u_1, v_1\}$
(seen as symbols rather than words).

Similarly, we can rewrite $w$ uniquely as a word over $\{u_k,v_k\}$ for all
$k$. We denote $\phi_k (w)$ this rewriting.

Consider now $N = \abs{u_{k+1}} + 4 \abs{u_k}$.
A factor $f$ of $w$ is right-special if it can be extended in two ways on
the right.
Given $f$ of length $N$, it can be rewritten as a word over the symbols
$\{u_k, v_k \}$, apart maybe from a prefix and a suffix, each of which is of
length at most $2 \abs{u_k}$ (the length is taken over $\{a,b\}$).
These prefix and suffix can possibly not be recognized, which means it
can't be determined whether they belong to $u_k$ or $v_k$.
Still, there is in $f$ at least one occurrence of one of the following patterns:
$(u_k)^4$ or $v_k (u_k)^3 v_k$.
The occurrence of either one of these patterns allows to recognize that $f$ is
a subfactor of $x u_{k+1}$ or $x v_{k+1}$, respectively, where
$x \in \{u_{k+1},v_{k+1}\}$.
Therefore in most cases, it allows us to predict how $f$ can be extended on the
right, and the extension is unique.
We know that $u_i$ and $v_i$ begin by $u_{i-1}$ for all $i$, and $u_1,v_1$ begin
by $a$.
So the only case when the extension is not unique is when $f$ contains
$x u_{k}u_{k-1} \ldots u_1 a$ as a right factor, with
$x \in \{u_{k+1},v_{k+1}\}$.
\begin{itemize}
 \item If the next letter is $a$, then the $aa$ at the end of $fa$ will
determine that the unique extension of $fa$ is
$x_{k+1}u_{k}u_{k-1} \ldots u_1 u_1$.
This will imply (by induction), that the unique extension of $fa$ is
$x_{k+1} u_{k+1}$.
 \item Similarly, if the next letter is $b$, the unique extension of $fb$ will
be $x_{k+1} v_{k+1}$.
\end{itemize}
We have to prove that a factor of length $N$ which ends by
$x_{k+1}u_{k}u_{k-1} \ldots u_1 a$ is unique: indeed, looking at the form of
$v_{k+1}$ and $u_{k+1}$, the finite word $f$ is a factor of
$v_k (u_k)^2 x_{k+1} u_{k} \ldots u_1 a$, which is of length greater
than $N$.

Therefore, there exists exactly one right-special word for this value of $N$.
\ep\medbreak

We remarked that in the case of a two-letter alphabet, the number of
right-special words was equal to $s(n)$.
Therefore, the lemma proves (thanks to proposition~\ref{bound-cohom-b}) that
infinitely often, the Rauzy graph of $w$ has an homology of rank two.
It means that it is homotopic to a wedge of two circles.
Application of proposition~\ref{bound-cohom-b} and lemma~\ref{lemma-cohom-ccl}
yields the following result:
\begin{prop}
 The suspension space associated with the word $w$ has finitely generated
cohomology over $\Q$.
\end{prop}

This last example proves that cut and projection tiling spaces, unlike other
tiling spaces, have some kind of rigidity which links complexity and
cohomology.

\section*{Acknowledgments}
I would like to thank in first place my advisor, Johannes Kellendonk, for
guidance and support.
I am also grateful to all the people with whom I could discuss about this work,
for their useful comments and advices: Boris Adamczewski, Nicolas Bedaride, Jean
Bellissard, Val\'erie Berth\'e, Julien Cassaigne, Jean-Marc Gambaudo, Laurent
Vuillon.

\end{document}